\documentclass{ert-l}%
\issueinfo{00}
  {}
  {}
  {1995}
 
\copyrightinfo{1995}
  {American Mathematical Society}
\usepackage{amssymb} 
\newcommand{\CC}{{\mathbb C}} 
\newcommand{\ZZ}{{\mathbb Z}}

\newcommand{\liea}[1]{\mathfrak{#1}} 
\newcommand{\lieg}[1]{\mathrm{#1}}

\newcommand{\ad}{\operatorname{ad}}

\theoremstyle{plain} 
\newtheorem{thm}{Theorem}[section]

\newtheorem{lm}[thm]{Lemma} 
 
\newtheorem{cor}[thm]{Corollary} 
 
\newtheorem{prop}[thm]{Proposition}

\theoremstyle{definition} 
\newtheorem{re}[thm]{Remark} 
\newtheorem{ex}[thm]{Example} 
 
\newtheorem{recipe}[thm]{Recipe} 
 
\begin{document} 
 
\title{Nice parabolic subalgebras of reductive Lie algebras} 
\author{Karin Baur and Nolan Wallach} 
\thanks{First named author supported by the Swiss National 
Science Foundation,\\ 
 Second named author partially suppoted by an NSF summer grant.} 
\address{Karin Baur, Department of Mathematics, University of 
California, San Diego, USA\\ 
Nolan Wallach, Department of Mathematics, University of 
California, San Diego, USA} 
\email{kbaur@math.ucsd.edu, nwallach@math.ucsd.edu} 
 
\date{October 5, 2004} 
 
\maketitle 
 
\section*{Abstract} 
This paper gives a classification of parabolic subalgebras of 
simple Lie algebras over $\CC$ that are complexifications of 
parabolic subalgebras of real forms for which Lynch's vanishing 
theorem for generalized Whittaker modules is non-vacuous. The 
paper also describes normal forms for the admissible characters in 
the sense of Lynch (at least in the quasi-split cases) and 
analyzes the important special case when the parabolic is defined 
by an even embedded TDS (three dimensional simple Lie algebra).

%
\section*{Introduction} 
%
%
If $\mathfrak{g}$ is a semi-simple Lie algebra over $\mathbb{C}$ 
and if $\mathfrak{p}$ is a parabolic subalgebra of $\mathfrak{g}$ 
then there is a $\mathbb{Z}$-grade of $\mathfrak{g}$ as a Lie 
algebra, $\mathfrak{g}= \sum_{j}\mathfrak{g}_{j}$ such that 
$\mathfrak{p}= \sum_{j\geq0}\mathfrak{g}_{j}$ and if 
$\mathfrak{n}= \sum_{j>0} \mathfrak{g}_{j}$\ then $\mathfrak{n}$ 
is the nilradical of $\mathfrak{p}$ and $\mathfrak{g}_{1}$ 
projects bijectively onto the abelianization 
($\mathfrak{n/[n},\mathfrak{n}]$) of $\mathfrak{n}$. The purpose 
of this article is to give a classification of those parabolic 
subalgebras such that there is a Richardson element 
$x\in\mathfrak{n}$ such that $x\in\mathfrak{g}_{1}$. That is there 
exists $x\in \mathfrak{g}_{1}$ such that 
$[\mathfrak{p},x]=\mathfrak{n}$. A parabolic subalgebra will be 
called \emph{nice} if it satisfies this condition. \ These 
parabolic subalgebras are exactly the complexifications of the 
real parabolic subalgebras whose nilradicals support admissible 
Lie algebra homomorphisms to $i\mathbb{R}$ in the sense of Lynch's 
thesis ~\cite{l}. In that thesis Lynch proved a generalization of 
Kostant's vanishing theorem for Whittaker modules [K] (valid for 
generic Lie algebra homomorphisms to $i\mathbb{R}$ of nilradicals 
of Borel subalgebras of quasisplit real forms).  Of course, Lynch's 
theorem is vacuous if the nilradical of the parabolic subalgebra 
admits no such homomorphisms. He introduced the term admissible 
for the parabolic sublagebras whose nilradicals  admit admissible 
homomorphisms. Thus an admissible parabolic in the sense of Lynch 
is a real form of a nice parabolic in our sense. 
 
It is clear that a parabolic subalgebra is nice if and only if its 
intersection with each ideal is nice. Thus it is 
enough to do the classification for simple Lie algebras over 
$\mathbb{C}$. The complete classification is given in section 1 
and the rest of the paper is devoted to the proof of the 
correctness of the list (sections 3 and 4), and to a 
description of 
the corresponding Richardson elements (admissible elements in the 
sense of Lynch). The proofs of the assertions in section 5 will 
appear in [B]. 
 
In ~\cite{wa} the second named author used his extension of the 
Lynch results to prove a holomorphic continuation of generalized 
Jacquet integrals for degenerate principal series under the real 
analogue of the condition of niceness. These results contain all 
known cases of continuation, holomorphy and (essentially) 
uniqueness for Jaquet integrals and Whittaker models. \ Thus the 
results of this paper explain both the range of applicability of 
those results and their limitations. 
 
In ~\cite{l} Lynch studies the classification of his admissible 
parabolic subalgebras. His results are correct for types 
$\lieg{A}_n$, $\lieg{F}_4$ and $\lieg{G}_2$. 
 
The paper [EK] studies the notion of a \emph{good grade} of a 
semi-simple algebra. \ This is a $\mathbb{Z}$-grade of 
$\mathfrak{g} =\sum_{j} \mathfrak{g}_{j}$ as a Lie algebra such 
that there is an element 
$x\in\mathfrak{g}_{2}$ such that $\ad(x)$ is injective on $%
\sum_{j<0} \mathfrak{g}_{j}$. \ Theorem ~\ref{th:surj} implies 
that a good grade with all odd components equal to $0$ defines a 
nice parabolic $ \sum_{j\geq0} \mathfrak{g}_{j}$. \ Thus a 
classification of good grades yields as a special case a 
classification of nice parabolic subalgebras. That said we have 
decided that the special case was of sufficient importance to have 
an independent exposition even if the paper [EK] were without any 
errors (see the remark before Theorem \ref{th:D}). \ Also, our 
description of the answer for the classical groups is (we hope) 
simpler than that given in [EK]. \ In addition, we have endeavored 
to give enough detail that a serious reader could with little 
additional effort check that the results are correct. 
 
The authors would like to thank Elashvili for informing them of 
his work with Kac. For obvious reasons there is substantial 
overlap between our papers. Many of the general results in section 
2 can be found in the second named author's manuscript [W] that 
has been freely available on his web site (in various forms) for 
at least two years and in [EK].

\section{Statement of the results}\label{sec:one} 
%
%
If not specified otherwise, $\liea{g}$ will denote a simple Lie 
algebra over the complex numbers. 
Fix a Borel subalgebra $\liea{b}$ in $\liea{g}$, let 
$\liea{h}\subset\liea{b}$ be a Cartan subalgebra of $\liea{g}$. We 
will denote the set of simple roots relative to this choice by 
$\Delta=\{\alpha_1,\dots,\alpha_n\}$. We always use the 
Bourbaki-numbering of simple roots. 

Let $\liea{p}\subset\liea{g}$ be a parabolic subalgebra, 
$\liea{p}=\liea{m}\oplus\liea{u}$ (where $\liea{m}$ is a Levi 
factor and $\liea{u}$ the corresponding nilpotent radical of 
$\liea{p}$). After conjugation we can assume that $\liea{p}$ 
contains the chosen Borel subalgebra and 
$\liea{m}\supset\liea{h}$. If $\liea{b}$ has been fixed then we 
will say that $\liea{p}$ is standard if  $\liea{p}\supset\liea{b}$ 
from now on. In particular, if $\liea{p}$ is standard then it is 
given by a subset of $\Delta$, namely the simple roots such that 
both roots spaces $\liea{g}_{\pm\alpha}$ belong to the Levi factor 
of $\liea{p}$. 
 
Thus such a parabolic subalgebra is described by a $n$-tuple, 
$(u_1,...,u_n)$ in $\{0,1\}^n$: ones correspond to simple roots 
with root spaces not in $\liea{m}$. Equivalently, a parabolic 
subalgebra is given by a coloring of the Dynkin diagram of the Lie 
algebra: a black (colored) node corresponds to a simple root whose 
root space belongs to $\liea{m}$. Here, one has to be very careful 
since there exist different notations. Our choice was motivated by 
the coloring for Satake diagrams.  Let $(u_1,...,u_n)$ define the 
parabolic subalgebra $\liea{p}$ and and $H \in \liea{h}$ be 
defined by $\alpha_i(H) = u_i$. If we set $\liea{g}_i=\{x\in 
\liea{g}|[H,x]=ix\}$ then $\liea{p}=\sum_{i\geq 0}\liea{g}_i$. 
 
%
%
\subsection{Results in the classical cases} 
%
%
As is usual, we will refer to the simple Lie algebras of type 
$\lieg{A}_n,\lieg{B}_n,\lieg{C}_n,\lieg{D}_n$ as the classical Lie 
algebras and the remaining five simple Lie algebras will be called 
exceptional.  We realize the classical Lie algebras as subalgebras 
of $\liea{gl}_N$ for $N=n+1, 2n+1,2n,2n$ respectively. With 
$\lieg{A}_n$ the trace zero matrices, $\lieg{B}_n,\lieg{D}_n$ the 
orthogonal Lie algebra of the symmetric form with matrix with all 
entries $0$ except for those on the skew diagonal which are $1$ 
and $\lieg{C}_n$ the symplectic Lie algebra for the symplectic 
form with matrix whose only nonzero entries are skew diagonal and 
the first $n$ are $1$ and the last $n$ are $-1$. With this 
realization we take as our choice of Borel subalgebra the 
intersection of the corresponding Lie algebra with the upper 
triangular matrices in $\liea{gl}_N$. We will call a parabolic 
subalgebra that contains this Borel subalgebra standard. If 
$\liea{p}$ is a standard parabolic subalgebra then we refer to the 
Levi factor that contains the diagonal Cartan subalgebra, 
$\liea{h}$, by $\liea{m}$ and call it the standard Levi factor. 
 
Thus for all classical Lie algebras the standard Levi factor is then 
in diagonal block form given by a sequence of square matrices on 
the diagonal. For the orthogonal and symplectic Lie algebras, 
these sequences are palindromic.  If $\liea{p}$ is a parabolic 
subalgebra for one of these Lie algebras then if $\liea{m}$ is the 
standard Levi factor of the parabolic subalgebra to which it is 
conjugate then we will say that $\liea{m}$ is the standard Levi 
factor. 

\begin{thm}\label{th:A} 
For type $\lieg{A}_n$ a parabolic subalgebra is nice if and only 
if the sequence of the block lengths of its standard Levi factor 
is unimodal. 
\end{thm} 
 
\begin{thm}\label{th:B} 
For type $\lieg{B}_n$ a parabolic subalgebra $\liea{p}$ is nice if 
and only if its standard Levi factor either has unimodal block 
lengths or the block lengths are of the form 
\[ a_1 \leq a_2 \leq ...\leq a_r > b_1=b_2=...=b_s < 
a_r \geq ...\geq a_1\] 
with $b_1 = a_r-1$ and if $r>1$ then $a_{r-1} < a_r$. 
\end{thm} 
 
\begin{thm}\label{th:C} 
For type $\lieg{C}_n$ a parabolic subalgebra is nice if and only if 
the sequence 
of block lengths of its standard Levi factor is unimodal and if 
the number of blocks 
is odd then each odd block length occurs exactly twice. 
\end{thm} 
 
Before we state our result for $\lieg{D}_n$ a word should be said 
about the ambiguity in describing parabolic subalgebras in the 
case of $\lieg{D}_n$. We note that the intersection of the 
standard parabolic subalgebra with block sizes 
$[a_1,a_2,...,a_r,1,1,a_r,...,a_2,a_1]$ for $\lieg{GL}(d)$ 
($d=2(a_1+...+a_r+1)$) with $\liea{so}(d)$ is the same as the 
intersection with the standard parabolic with block lengths 
$[a_1,a_2,...,a_r,2,a_r,...,a_2,a_1]$. In the following theorem we 
we will only look at the second version of the parabolic subgroup. 
The first two types in Theorem 6.4 (i),(ii) in [EK] do not appear 
in our result directly because of this choice. 
 
We should also point out that there is a difference between our 
statement and that of [EK] involving the last part of this result 
and Theorem 6.4 (i) and (ii) in~\cite{ek}. Their list is missing 
all cases with $s$ even and $a_r > 3$.  In addition to having a 
(perhaps too) detailed proof of the following result in this paper 
the authors have done extensive computer computations which agree 
with our formulation for $n \leq 20$. 
\begin{thm}\label{th:D} 
For type $\lieg{D}_n$ a parabolic subalgebra is nice if and only 
if its standard Levi factor (taking into account the choice made 
above) has one of the following forms: 
 
1) It is  unimodal with an odd number of blocks. 
 
2) It is unimodal with an even number of block lengths and the odd 
block lengths occur exactly twice. 
 
3) The block lengths are of the form 
\[ a_1 \leq a_2 \leq ...\leq a_r > b_1=b_2=...=b_s < a_r \geq ...\geq a_1\] 
with $b_1 = a_r-1$, $a_r$ is odd and $a_{r-1}<a_r$ 
and if $s$ is even then the 
odd block lengths occur exactly twice. 
\end{thm} 
 
\subsection{Results in the exceptional cases} 
In this subsection we will state the classification of 
nice parabolic subalgebras for the 
exceptional simple Lie algebras. The parabolic subalgebra will be 
given by an $n$-tuple where $n$ is the rank and the entries are 
$\alpha_i(H)$ where $H$ is the element that gives the grade 
corresponding to the parabolic subalgebra and the $\alpha_i$ are 
the simple roots in the Bourbaki order. 
 
The only nice parabolic subalgebras in type $\lieg{G}_2$, 
$\lieg{F}_4$ are those that are given by an even 
$\liea{sl}_2$-triple (see the next section for the definition). 
They are listed below. 
\[ 
\begin{array}{cc} 
\lieg{G}_2 & \lieg{F}_4 \\ 
 & \\ 
(1,1) & (1,1,1,1)\\ 
(1,0) & (1,1,0,1)\\ 
(0,0) & (1,1,0,0) \\ 
& (1,0,0,1) \\ 
& (0,1,0,1) \\ 
& (0,1,0,0) \\ 
& (0,0,0,1) \\ 
& (0,0,0,0) 
\end{array} 
\] 
The nice parabolic subalgebras of type $\lieg{E}$ are given in the 
following table: 
\[ 
\begin{array}{rccc} 
 & \lieg{E}_6 & \lieg{E}_7 & \lieg{E}_8 \\ 
 & & & \\ 
1\ &(1,1,1,1,1,1) & (1,1,1,1,1,1,1) & (1,1,1,1,1,1,1,1) \\ 
2\ &(1, 1, 1, 0, 1, 1) & (1, 1, 1, 0, 1, 1, 1) 
& (1, 1, 1, 0, 1, 1, 1, 1) \\ 
3\ &(1, 1, 1, 0, 1, 0) & (1, 1, 1, 0, 1, 0, 1) 
& (1, 1, 1, 0, 1, 0, 1, 1) \\ 
4\ &(1, 1, 0, 1, 0, 1) & (1, 1, 0, 0, 1, 0, 1) 
&  \\ 
5\ &(1, 1, 0, 0, 1, 0) & (1, 1, 0, 0, 0, 0, 1) 
& (1, 0, 0, 1, 0, 1, 1, 1) \\ 
6\ &(1, 1, 0, 0, 0, 1) & (1, 0, 1, 1, 0, 1, 0) 
& (1, 0, 0, 1, 0, 1, 0, 1) \\ 
7\ &(1, 1, 0, 0, 0, 0) & (1, 0, 1, 0, 0, 1, 0) 
& (1, 0, 0, 1, 0, 0, 1, 1) \\ 
8\ &(1, 0, 1, 1, 0, 1) & (1, 0, 1, 0, 0, 0, 0) 
& (1, 0, 0, 1, 0, 0, 1, 0) \\ 
9\ &(1, 0, 1, 0, 0, 1) & (1, 0, 0, 1, 0, 1, 1) 
& (1, 0, 0, 0, 1, 0, 0, 1) \\ 
10\ &(1, 0, 1, 0, 0, 0) & (1, 0, 0, 1, 0, 1, 0) 
& (1, 0, 0, 0, 0, 1, 1, 1) \\ 
11\ &(1, 0, 0, 1, 1, 1) & (1, 0, 0, 1, 0, 0, 1) 
& (1, 0, 0, 0, 0, 1, 0, 1) \\ 
12\ &(1, 0, 0, 1, 0, 1) & (1, 0, 0, 0, 1, 0, 0) 
& (1, 0, 0, 0, 0, 1, 0, 0) \\ 
13\ &(1, 0, 0, 1, 0, 0) & (1, 0, 0, 0, 0, 1, 1 ) 
& (1, 0, 0, 0, 0, 0, 1, 1) \\ 
14\ &(1, 0, 0, 0, 1, 1) & (1, 0, 0, 0, 0, 1, 0) 
& (1, 0, 0, 0, 0, 0, 1, 0) \\ 
15\ &(1, 0, 0, 0, 1, 0) & (1, 0, 0, 0, 0, 0, 1) 
& (1, 0, 0, 0, 0, 0, 0, 1) \\ 
16\ &(1, 0, 0, 0, 0, 1) & (1, 0, 0, 0, 0, 0, 0) 
& (1, 0, 0, 0, 0, 0, 0, 0) \\ 
17\ &(1, 0, 0, 0, 0, 0) & (0, 1, 1, 0, 0, 1, 1) 
& (0, 1, 0, 0, 0, 0, 0, 1) \\ 
18\ &(0, 1, 1, 0, 1, 1) & (0, 1, 0, 0, 0, 0, 0) 
& (0, 1, 0, 0, 0, 0, 0, 0) \\ 
19\ &(0, 1, 1, 0, 0, 1) & (0, 0, 1, 0, 0, 1, 0) 
& (0, 0, 1, 0, 0, 0, 1, 0) \\ 
20\ &(0, 1, 0, 1, 0, 0) & (0, 0, 1, 0, 0, 0, 1) 
& (0, 0, 0, 1, 0, 0, 1, 1) \\ 
21\ &(0, 1, 0, 0, 0, 1) & (0, 0, 1, 0, 0, 0, 0) 
& (0, 0, 0, 1, 0, 0, 1, 0) \\ 
22\ &(0, 1, 0, 0, 0, 0) & (0, 0, 0, 1, 0, 1, 0) 
& (0, 0, 0, 1, 0, 0, 0, 1) \\ 
23\ &(0, 0, 1, 0, 0, 1) & (0, 0, 0, 1, 0, 0, 1) 
& (0, 0, 0, 0, 1, 0, 0, 1) \\ 
24\ &(0, 0, 1, 0, 0, 0) & (0, 0, 0, 1, 0, 0, 0) 
& (0, 0, 0, 0, 1, 0, 0, 0) \\ 
25\ &(0, 0, 0, 1, 0, 1) & (0, 0, 0, 0, 1, 0, 1) 
& (0, 0, 0, 0, 0, 1, 0, 0) \\ 
26\ &(0, 0, 0, 1, 0, 0) & (0, 0, 0, 0, 1, 0, 0) 
& (0, 0, 0, 0, 0, 0, 1, 1) \\ 
27\ &(0, 0, 0, 0, 1, 1) & (0, 0, 0, 0, 0, 1, 0) 
& (0, 0, 0, 0, 0, 0, 1, 0) \\ 
28\ &(0, 0, 0, 0, 1, 0) & (0, 0, 0, 0, 0, 0, 1) 
& (0, 0, 0, 0, 0, 0, 0, 1) \\ 
29\ &(0, 0, 0, 0, 0, 1) & (0, 0, 0, 0, 0, 0, 0) 
& (0, 0, 0, 0, 0, 0, 0,0) \\ 
30\ &(0, 0, 0, 0, 0, 0) &  & 
\end{array}\]

%
%
%
 
%
\section{Characterizations of niceness} 
%
%
In this section we study the notion of niceness and prove 
properties of nice parabolic subalgebras. Recall that a parabolic 
subalgebra induces a $\ZZ$-grading of $\liea{g}$ via the map $\ad 
H$ where $H\in\liea{h}$ is defined as in the beginning of section 
1.  We will also use the notation $B$ for the Killing form of 
$\liea{g}$. 
\begin{thm}\label{th:surj} 
Let $\liea{p}\subset\liea{g}$ with associated grading 
$\liea{g}=\oplus_{j\in\ZZ}\liea{g}_j$. The following are 
equivalent: 
\begin{enumerate} 
\item $\liea{p}$ is nice \item There exists $X\in\liea{g}_1$ such 
that $\ad(X):\liea{g}_j\to\liea{g}_{j+1}$ is surjective for all $j 
\geq 0$. \item There exists $X\in\liea{g}_1$ such that 
$\ad(X):\liea{g}_j\to\liea{g}_{j+1}$ is surjective for all $j > 
0$. 
\item There exists $X\in\liea{g}_1$ such that 
$\ad(X):\liea{g}_{j-1}\to\liea{g}_j$ is injective for all $j\leq 
0$. 
\end{enumerate} 
\end{thm} 
\begin{proof} 
Let $X\in\mathfrak{g}_{1}$ be such that $[\mathfrak{p},X]=\mathfrak{n}$. Then 
since \[\mathfrak{p=g}_{0}\oplus\mathfrak{g}_{1}\oplus...\oplus\mathfrak{g}%
_{r}\] we must have $\ad(X)\mathfrak{g}_{i}=\mathfrak{g}_{i+1}$ 
for $i\geq0$. This proves the necessity of condition 2). The 
sufficiency is equally clear. To prove 3) we will show that 3) 
implies 2). Let $\Omega$ be the set of all $X \in \mathfrak{g}_1$ 
satisfying the condition in 3). Then $\Omega$ is Zariski open and 
dense in $\mathfrak{g}_1$. We also note that if $\Lambda$ is the 
set of all $X \in \mathfrak{g}_1$ such that 
$[X,\mathfrak{p}]=\mathfrak{g}_1$ then it is also Zariski open and 
dense.  Hence $\Omega \cap \Lambda \neq \emptyset$. We have proved 
that 3) implies 2). 
 
We now prove 4). For this we observe that if 
$\overline{\mathfrak{n}}= 
{\textstyle\sum\limits_{i<0}}\mathfrak{g}_{i}$ (resp. 
$\overline{\mathfrak{p}}={\textstyle\sum\limits_{i\leq 
0}\mathfrak{g}_{i}}$) then $y\in\overline{\mathfrak{n}}$ (resp. 
$y\in \overline{\mathfrak{p}}$) is $0$ if and only if 
$B(y,\mathfrak{n})=\{0\}$ (resp. $B(y,\mathfrak{p})=\{0\}$). Let 
$X\in\mathfrak{g}_{1}$ be such that 
$[\mathfrak{p},X]=\mathfrak{n.}$ Then if $y\in\overline 
{\mathfrak{n}}$ and $[X,y]=0$ then 
\[ 
\{0\}=B([X,y],\mathfrak{p})=B(y,[X,\mathfrak{p}])=B(y,\mathfrak{n}). 
\] 
So $y=0$. Suppose that $X\in\mathfrak{g}_{1}$ and $\ad{X}:\overline{\mathfrak{n}%
}\rightarrow\mathfrak{g}$ is injective then we assert that $\ad{X}:\mathfrak{p}%
\rightarrow\mathfrak{n}$ is surjective. Indeed, if not there would exist 
$y\in\overline{\mathfrak{n}}$ with $B(y,[\mathfrak{p},X])=\{0\}$ but $y\neq0$. 
But  $B(y,[\mathfrak{p},X])=B([X,y],\mathfrak{p}$) and since $[X,\overline 
{\mathfrak{n}}]\subset\overline{\mathfrak{p}}$ we would have the contradiction 
$y\neq0$ but $[X,y]=0$. 
\end{proof} 
As a consequence we have the following useful criterion for the 
exceptional Lie algebras. 
 
\begin{cor}\label{cor:dim} 
If $\liea{p}\subset\liea{g}$ is nice then 
 
(i) $\dim\liea{g}_j\ge\dim\liea{g}_{j+1}$ for 
all $j> 1$. 
 
(ii) $\dim\liea{g}_1>\dim\liea{g}_2$. 
\end{cor} 
\begin{proof} 
Part (2) of Theorem~\ref{th:surj} gives (i). 
Note that the line through $X$ is in the kernel 
of the map 
$\ad X:\liea{g}_1\to\liea{g}_2$. This gives part (ii). 
\end{proof} 
Any nonzero nilpotent element $x$ of a complex semisimple Lie 
algebra is the nilpositive element of an $\liea{sl}_2$-triple 
$\{x,h,y\}$ (Jacobson-Morozov Theorem, cf. Section 3.3 
in~\cite{cm}). The next result describes a large class of examples 
of nice parabolic subalgebras (most of them in many simple Lie 
algebras) . 
 
\begin{lm}\label{lm:even} 
Let $\liea{p}\subset\liea{g}$ with grading 
$\liea{g}=\oplus_{j\in\ZZ}\liea{g}_j$. Assume that there is a 
nonzero $x\in\liea{g}_1$ and $y\in\liea{g}$ such that the 
$\liea{sl}_2$-triple of $x$ is $\{x,2H,y\}$. Then $\liea{p}$ is 
nice. 
\end{lm} 
\begin{proof} 
This follows from part (4) of Theorem~\ref{th:surj} since the 
kernel of $\ad{x}$ is contained in the sum of the eigenspaces for 
the eigenvalues of $2H$ that are non-negative integers. 
\end{proof} 
If $\liea{p}\subset\liea{g}$ is a subalgebra satisfying the 
assumptions of Lemma~\ref{lm:even}, we will say that $\liea{p}$ is 
given by an even TDS (by a three-dimensional simple subalgebra). 
 
We recall a helpful result of~\cite{wa}: 
Let $\liea{g}$ be a simple Lie algebra with Dynkin 
diagram $\triangle$. Let $\triangle'$ be a connected 
subdiagram of $\triangle$ and $\liea{g}'$ the corresponding 
simple Lie algebra. If $\liea{p}\subset\liea{g}$ is 
a parabolic subalgebra we denote the parabolic 
subalgebra $\liea{p}\cap\liea{g}'$ by $\liea{p}'$. 
 
\begin{lm}\label{lm:subdiag} 
If $\liea{p}\subset\liea{g}$ is nice then for 
every $\liea{g}'$ obtained as above, 
$\liea{p'}$ is nice. 
\end{lm} 
 
\begin{proof} 
Let $\liea{n}'$ be the nilradical of $\liea{p}'$ and let 
$\liea{n}$ be the nilradical of $\liea{p}$. Let $u \in \liea{h}$ 
be such that if $\alpha \notin \Delta'$ then $\alpha(u) = 1$. Let 
$\liea{v}$ be the direct sum of the eigenspaces for strictly 
positive eigenvalues of $\ad{u}$. Then $\liea{v}$ is an ideal in 
$\liea{n}$ and since $\liea{n}=\liea{n}' \oplus \liea{v}$, 
$\liea{n}' \cong \liea{n}/\liea{v}$. Let $p$ denote the projection 
from $\liea{n}$ onto $\liea{n}/\liea{v}$.  Let $X \in \liea{g}_1$ 
be as in (2) of Theorem~\ref{th:surj}.  Let $X'$ be the unique 
element in $\liea{n}'$ such that $p(X')=p(X)$. Then it is easy to 
see that $X' \in \liea{g}'_1$ and satisfies (3) of 
Theorem~\ref{th:surj}. 
\end{proof} 
 
Lemma~\ref{lm:subdiag} has two consequences: on one hand any 
subdiagram of a nice parabolic subalgebra is again nice. On the 
other hand, if a diagram contains a subdiagram that is not nice, 
then the ``big'' diagram is not nice. 
 
The following result will be one of our main tools in our proofs 
of the theorems in the previous section. 
\begin{thm}\label{thm:centr} 
Let $\liea{p}$ be a parabolic subalgebra of $\liea{g}$, 
$\liea{p}=\liea{m}\oplus\liea{u}$ where $\liea{m}$ is a Levi 
factor of $\liea{p}$ and $\liea{u}$ the corresponding nilpotent 
radical. Let $\liea{g} = \sum \liea{g}_j$ be the corresponding 
grade. Then 
 
$x\in\liea{g}_1$ is a Richardson element of $\liea{p}$ if and only 
if $\dim\liea{g}^x=\dim\liea{m}$. 
\end{thm} 
\begin{proof} 
We denote the opposite nilradical, $\sum_{j<0} \liea{g}_j$, by 
$\overline{\liea{u}}$. If $x\in\liea{u}$ then 
$\ad(x)\liea{g}=\ad(x)\overline{\liea{u}}+\ad(x)\liea{p}$. Now 
$\ad(x)\liea{p}\subset\liea{u}$ and 
$\dim\ad(x)\overline{\liea{u}}\le\dim\overline{\liea{u}}$. Thus 
\[ 
\dim\ad(x)\liea{g}\le\,2\dim\liea{u}. 
\] 
This implies for $x\in\liea{u}$ that 
$\dim\liea{g}^x\ge\dim\liea{m}$ and equality implies that 
$\dim\ad(x)\liea{p}=\dim\liea{u}$. Thus equality implies that $x$ 
is a Richardson element. If $x \in \liea{g}_1$ then 
$\ad(x)\overline{\liea{u}}\subset \sum_{j \leq 0} \liea{g}_j$. If 
$x$ is a Richardson element then since $\ad(x)$ is injective on 
$\overline{\liea{u}}$ the kernel of $\ad(x)$ is contained in 
$\liea{p}$ (see Theorem ~\ref{th:surj} (4)). Since 
$\ad(x)\liea{p}=\liea{u}$ we have proved that $\dim \liea{g}^x = 
\dim \liea{m}$. 
\end{proof} 
 
We recall that {\it distinguished} nilpotent elements 
of $\liea{g}$ are nilpotent elements that 
are not contained in any Levi factor 
$\liea{m}\subsetneq\liea{g}$. 
Let $x\in\liea{g}_1$ 
be nilpotent, let $x,2H,y$ be a TDS and 
$\liea{g}=\sum\liea{g}_i$ the grading given by $H$. 
Then $x$ 
is distinguished if and only if 
$\dim\liea{m}=\dim\liea{g}_1$ (cf. Section 8.2 of~\cite{cm}). 
 
\begin{cor}\label{cor:dist} 
Let $\liea{p}\subsetneq\liea{g}$ and $x\in\liea{g}_1$ 
a Richardson element with $\dim\liea{m}>\dim\liea{g}_1$. 
Then $x$ is not a 
distinguished nilpotent element for $\liea{g}$. 
\end{cor} 
 
\begin{proof} 
Follows directly from Theorem~\ref{thm:centr}. 
\end{proof} 
 
\section{The determination of Jordan forms}\label{sec:jordan} 
 
We use the notation $M_{p,q}(\mathbb{C})$ for the $p\times q$ 
matrices over $\mathbb{C}$. Let $\ n_{1},...,n_{r}>0$ with 
$r\geq2$ be given. Then we define a map 
$\phi:M_{n_{1},n_{2}}(\mathbb{C})\times M_{n_{2},n_{3}}(\mathbb{C}%
)\times\cdots\times M_{n_{r-1},n_{r}}(\mathbb{C})\rightarrow M_{n_{1},n_{r}%
}(\mathbb{C}$) by $\phi(X_{1},...,X_{r-1})=X_{1}\cdots X_{r-1} $. 
The following should be well known. Since it will be used several 
times we include a proof. 
 
\begin{lm}\label {lm:basic} 
If $r\geq2$ then $\phi$ is a surjection onto the variety $Y_{n_{1}%
,n_{r},m}=\{X\in M_{n_{1},n_{r}}|$rank$(X)\leq m\}$ with $m=\min 
\{n_{1},...,n_{r}\}.$ 
\end{lm} 
\begin{proof} 
If $V$ is a vector space over $\mathbb{C}$ then we will use the 
notation $\mathbb{P}V$ for the projective space of $V$. We observe 
that $\phi$ induces a 
regular map of $\mathbb{P}M_{n_{1},n_{2}}(\mathbb{C})\times\mathbb{P}%
M_{n_{2},n_{2}}(\mathbb{C})\times\cdots\times\mathbb{P}M_{n_{r-1},n_{r}%
}(\mathbb{C})$ into $\mathbb{P}M_{n_{1},n_{r}}(\mathbb{C})$. We 
will prove that the image of this map is the projective variety 
$Z=\{[X]\in \mathbb{P}M_{n_{1},n_{r}}(\mathbb{C})|\Delta(X)=0$ for 
all $m+1\times m+1$ minors of $X\}$. This will clearly prove the 
Lemma. \ We first observe that $Y=Y_{n_{1},n_{r},m}$ is an 
irreducible affine variety. Indeed, define the $n_{1}\times n_{r}$ 
matrix $J$ by 
\[ 
J=\left[ 
\begin{array} 
[c]{cc}%
I & 0\\ 
0 & 0 
\end{array} 
\right]  , 
\] 
with $I$ the $m\times m$ identity matrix. Then the set of matrices 
in $M_{n_{1},n_{r}}(\mathbb{C})$ of rank exactly $m$ is 
$U=\{gJh^{-1}|g\in GL(n_{1},\mathbb{C}),h\in 
GL(n_{r},\mathbb{C})\}$. Since $U$ is irreducible as a 
quasi-projective variety (indeed it is affine) and $Y$ is the 
Zariski closure of $U$ the irreducibility follows. 
 
We now prove the result by induction on $r$. If $r=2$ then $\phi$ 
is the identity map and the result is obvious. Assume for $r$ and 
we will now prove it for $r+1$. \ By the above and the inductive 
hypothesis it is enough tp show 
that the set $H=\{WX|X\in M_{n_{r},n_{r+1}}(\mathbb{C}),W\in Y_{n_{1},n_{r}%
,m}$ is Zariski dense in $Y_{n_{1},n_{r+1},m^{\prime}}$ with $m^{\prime}%
=\min\{n_{r+1,}m\}$. To see this we observe that if $J$ is as 
above and if $K$ is $\left[ 
\begin{array} 
[c]{cc}%
I & 0 
\end{array} 
\right]  $ if $n_{r}\leq n_{r+1}$ or $\left[ 
\begin{array} 
[c]{c}%
I\\ 
0 
\end{array} 
\right]  $ if $n_{r}>n_{r+1}$ then the set $\{gJKh^{-1}|g\in GL(n_{1}%
,\mathbb{C}),h\in GL(n_{r+1},\mathbb{C})\}$ is contained in $H$ 
and $JK$ is the \textquotedblleft$J$\textquotedblright\ for 
$Y_{n_{1},n_{r+1},m^{\prime}}$. 
\end{proof} 
We now consider a standard parabolic subalgebra $\mathfrak{p}$ of 
A$_{n}$, $n\geq1$. Then $\mathfrak{p}$ is given by integers 
$n_{1},...,n_{r}>0$ with $n_{1}+...+n_{r}=n+1$. We think of the 
matrices as given in $r\times r$ block form with the $i$-th 
diagonal block of size $n_{i}\times n_{i}$. \ Then the 
grade corresponding to $\mathfrak{p}$ is given by bracket with%
\[ 
H=\left[ 
\begin{array} 
[c]{cccc}%
\frac{r-1}{2} I_{n_1} & 0 & 0 & 0\\ 
0 & \frac{r-3}{2} I_{n_2}& 0 & 0\\ 
\vdots & \vdots & \ddots & \vdots\\ 
0 & 0 & 0 & -\frac{r-1}{2} I_{n_r}%
\end{array} 
\right]  . 
\] 
Thus $x\in\mathfrak{g}_{1}$ is of the block form%
\[ 
x=\left[ 
\begin{array} 
[c]{ccccc}%
0 & X_{1} & 0 & 0 & 0\\ 
0 & 0 & X_{2} & 0 & 0\\ 
\vdots & \vdots & \vdots & \ddots & \vdots\\ 
0 & 0 & 0 & 0 & X_{r-1}\\ 
0 & 0 & 0 & 0 & 0 
\end{array} 
\right] 
\] 
with $X_{i}$ a $n_{i}\times n_{i+1}$ matrix. If we consider the 
$j$-th power of $x$ then it has all of its (possibly) non-zero 
entries in the $j+1$-st super-diagonal and those entries are 
$X_{1}X_{2}\cdots X_{j},X_{2}X_{3}\cdots 
X_{j+1},...,X_{r-j}X_{r+1-j}\cdots X_{r-1}$. We have 
\begin{lm}\label{lm:rankA} 
If $x\in\mathfrak{g}_{1}$ is generic then $x^{j}$ has rank equal 
to 
\[%
{\textstyle\sum\limits_{i=1}^{r-j+1}} 
\min\{n_{i},n_{i+1},...,n_{i+j-1}\}. 
\] 
 
\end{lm} 
This is a simple application of the above observation, and Lemma 
~\ref{lm:basic}. 
 
We now recall an algorithm for calculating the Jordan canonical 
form of a nilpotent matrix. \ Let $x$ be a nilpotent $n+1\times 
n+1$ matrix. We assume that its Jordan canonical form has $a_{s}$ 
blocks of size $s$ for $s=1,...,m$ and $x^{m}=0$. If $J$ is a 
nilpotent Jordan block of size $k$ then the dimension of the 
kernel of $J^{s}$ is $s$ for $s=1,...,k-1$ and $k$ for $s\geq 
k$. Hence $a_{1},...,a_{r}$ satsify the system of linear equations%
\[ 
a_{1}+a_{2}+a_{3}+...+a_{r-1}+a_{r}=\dim\ker x 
\]%
\begin{align*} 
a_{1}+2a_{2}+2a_{3}+...+2a_{r-1}+2a_{r}  &  =\dim\ker x^{2}\\ 
&  \vdots 
\end{align*}%
\[ 
a_{1}+2a_{2}+3a_{3}+...+(r-1)a_{r-1}+(r-1)a_{r}=\dim\ker x^{r-1}%
\]%
\[ 
a_{1}+2a_{2}+3a_{3}+...+(r-1)a_{r-1}+ra_{r}=\dim\ker x^{r}=n+1. 
\] 
This system has a unique solution%
\[ 
a_{j}=-\dim\ker x^{j+1}+2\dim\ker x^{j}-\dim\ker x^{j-1}\underset{}%
{\qquad(\ast)}%
\] 
for $j=1,2,...,r$. 
 
If $x$ is a generic element in $\mathfrak{g}_{1}$ then the 
previous lemma 
implies that%
\[ 
\dim\ker x^{j}=n+1-%
{\textstyle\sum\limits_{i=1}^{r-j}} 
\min\{n_{i},n_{i+1},...,n_{i+j}\}. 
\] 
We have proved 
\begin{prop}\label{prop:jordA} 
Let $\mathfrak{p}$ be the parabolic subalgebra of A$_{n}$, 
$n\geq1$ given by integers $n_{1},...,n_{r}>0$ with 
$n_{1}+...+n_{r}=n+1$. Let $x$ be a generic element of 
$\mathfrak{g}_{1}$ (i.e. $[\mathfrak{p},x]=\mathfrak{g}_{1}$) then 
the Jordan canonical form of $x$ has $a_{j}$ blocks of size $j$ 
for $j=1,...,r$ with 
\begin{align*} 
a_{j}  &  =%
{\textstyle\sum\limits_{i=1}^{r-j-1}} 
\min\{n_{i},n_{i+1},...,n_{i+j+1}\}-2%
{\textstyle\sum\limits_{i=1}^{r-j}} 
\min\{n_{i},n_{i+1},...,n_{i+j}\}+\\ 
& 
{\textstyle\sum\limits_{i=1}^{r-j+1}} 
\min\{n_{i},n_{i+1},...,n_{i+j-1}\}. 
\end{align*} 
 
\end{prop} 
 
We will now study the Jordan forms for the other classical Lie 
algebras. We will use the versions of these Lie algebras as they 
are given in ~\cite{gw}. We will now give the analogue of 
Proposition \ref{prop:jordA} for these cases. We start with the 
case of $\lieg{C}_n$. 
 
\begin{prop}\label{prop:rank-C} 
Let $\mathfrak{g}$ be of type $\lieg{C}_{n}$ and let 
$\mathfrak{p}$ is a standard parabolic subalgebra. Let 
$x\in\mathfrak{g}_{1}$ be a generic element.  If $\mathfrak{p}$ is 
given by an even number of blocks then the element $x$ is generic 
for the corresponding parabolic subgroup of $\mathfrak{sl}_{2n}$ 
and thus its Jordan form is calculated as in Proposition 
\ref{prop:jordA}. If there are an odd number of blocks the Jordan 
form is calculated as follows: write the block sizes as 
$n_{1},...,n_{r},n_{r+1},n_{r+2},...,n_{2r+1}$ and denote \ this 
$2r+1$ tuple by $\mathbf{n}$. We define\ a non-negative integer $r_{ij}%
(\mathbf{n})$ as follows. If $i<1$, $j<1$ or $i+j>2r+1$ then $r_{ij}%
(\mathbf{n})=0$. We now assume $i,j\geq1$ and $i+j\leq2r+1$. If 
$i+j\leq r+1$ or if $i\geq r+1$ then 
$r_{ij}(\mathbf{n})=\min\{n_{i},n_{i+1},...,n_{i+j}\}$. If $i\leq 
r$ and $i+j>r+1$ then we set 
$s=\min\{n_{i},n_{i+1},...,n_{i+j}\}$, 
$u=\min\{r+1-i,i+j-r-1\}$ \ and $v=\min\{n_{r+1-u},n_{r+2-u},...,n_{r+1+u}%
\}$.  If $v$ is odd set $t=v-1$ if $v$ is even set $t=v$, then $r_{ij}%
(\mathbf{n)}=\min\{s,t\}$. Then the rank of $x^{j}$ is 
\[ 
r_{j}(\mathbf{n)=}\sum_{i}r_{ij}(\mathbf{n}),j>0,r_{0}(\mathbf{n)}=2n. 
\] 
For each $1\leq j\leq2r+1$ there are $r_{j-1}(\mathbf{n})-2r_{j}%
(\mathbf{n)+}r_{j+1}(\mathbf{n)}$ blocks of size $j$ in the Jordan 
form of $x$. 
\end{prop} 
 
The result for $\lieg{B}_{n}$ and $\lieg{D}_{n}$ is similarly 
complicated.  Also we will follow the convention introduced in the 
discussion before Theorem \ref{th:D} to resolve the ambiguity in 
the cases where there are an even number of block sizes and the 
central pair are both of size 1. These cases are ignored and 
we use the description of the same parabolic now with one less 
block size and all sizes the same except that the 1,1 is replaced 
by a 2. 
 
\begin{prop}\label{prop:rank-BD} 
Let $\mathfrak{g}$ be of type $\lieg{B}_{n}$ or $\lieg{D}_{n}$ and 
let $\mathfrak{p}$ is a standard parabolic subalgebra and let 
$d=2n+1$ in the case of $B_{n}$ and $d=2n$ for $D_{n}$. Let 
$x\in\mathfrak{g}_{1}$ be a generic element. If $\mathfrak{p}$ is 
given by an odd number of blocks then the element $x$ is generic 
for the corresponding parabolic subgroup of $\mathfrak{sl}_{d}$ and 
thus its Jordan form is calculated as in Proposition 
\ref{prop:jordA}. If there are an even number of blocks the Jordan 
form is calculated as follows: write the block sizes as 
$n_{1},...,n_{r},n_{r+2},...,n_{2r+1}$ and denote \ this $2r$ 
tuple by $\mathbf{n}$. We define\ a non-negative integer 
$r_{ij}(\mathbf{n})$ as follows. If $i<1$, $j<1$ or $i+j>2r$ then 
$r_{ij}(\mathbf{n})=0$. We now assume $i,j\geq1$ and $i+j\leq2r$. 
If $i+j\leq r$ or if $i\geq r$ then 
$r_{ij}(\mathbf{n})=\min\{n_{i},n_{i+1},...,n_{i+j}\}$. If $i\leq 
r$ and $i+j\geq r+1$ then we set 
$s=\min\{n_{i},n_{i+1},...,n_{i+j}\}$, $u=\min\{r+1-i,i+j-r\}$ \ 
and $v=\min\{n_{r+1-u},n_{r+2-u},...,n_{r+u}\}$. If 
$v$ is odd set $t=v-1$ if $v$ is even set $t=v$, then $r_{ij}(\mathbf{n)}%
=\min\{s,t\}$. Then the rank of $x^{j}$ is 
\[ 
r_{j}(\mathbf{n)=}\sum_{i}r_{ij}(\mathbf{n}),j>0,r_{0}(\mathbf{n)}=2n. 
\] 
For each $1\leq j\leq2r+1$ there are $r_{j-1}(\mathbf{n})-2r_{j}%
(\mathbf{n)+}r_{j+1}(\mathbf{n)}$ blocks of size $j$ in the Jordan 
form of $x$. 
\end{prop}

The intricacy of these statements is evidenced in the complicated 
proof. We warn the reader that it will be painful. We will devote 
the rest of this section to the proofs of these propositions. For 
B$_{n}$ and D$_{n}$ we will use the form with matrix having all 
entries $0$ except for ones on the skew diagonal%
\[ 
L=\left[ 
\begin{array} 
[c]{cccc}%
0 & 0 & 0 & 1\\ 
0 & 0 & 1 & 0\\ 
\vdots & \ddots & \vdots & \vdots\\ 
1 & 0 & 0 & 0 
\end{array} 
\right] 
\] 
For C$_{n}$ we use the form given by%
\[ 
M=\left[ 
\begin{array} 
[c]{cc}%
0 & L\\ 
-L & 0 
\end{array} 
\right] 
\] 
with $L$ as above. We will use the notation $A^{\#}$ for the 
adjoint of the matrix $A$ with respect to the indicated form. \ 
Thus for $\lieg{B}_{n},\lieg{D}_{n}$ if $A $ has $ij$ entry 
$a_{ij}$ then $A^{\#}$ has $ij$ entry $a_{m-j,m-i}$ with $m=2n+2$, 
$2n+1$ respectively. In the case $\lieg{C}_{n}$ we write 
$\eta_{i}=1$ if $i\leq n$ and $\eta_{i}=-1$ if $i>n$. \ Then if 
$A$ has $ij$ entry $a_{ij}$ then $A^{\#}$ has $ij$ entry 
$-\eta_{i}\eta_{j}a_{2n+1-j,2n+1-i}$. In all cases, if 
$e_{i},i=1,...,q$ is the standard basis with $q=2n+1,2n,2n$ for 
$\lieg{B}_{n},\lieg{C}_{n}$ and $\lieg{D}_{n}$ respectively than 
if $(...,...)$ is the 
corresponding invariant form on $\mathbb{C}^{q}$ then $(e_{i},e_{q+1-j}%
)=\pm\delta_{ij}$. This implies the following observations: 
 
\smallskip 
 
\noindent(a) If $i+p\leq n+1$ then $(...,...)$ induces a perfect 
pairing between 
\[V_{i,p}= 
\lieg{span}_{\mathbb{C}}\{e_{i},e_{i+1},...,e_{i+p-1}\}\] and 
\[\check{V}_{i,p} = \lieg{span}_{\mathbb{C}}\{e_{q+2-i-p},e_{q+1-i-p},...,e_{q+1-i}%
\}\] 
 
\smallskip 
 
\noindent(b) In the case of $\lieg{B}_{n}$ the form $(...,...)$ 
induces a nondegenerate form on 
span$_{\mathbb{C}}\{e_{i}|n+1-p\leq i\leq n+1+p\}$ for 
$p=0,...,n$. 
 
\smallskip 
 
\noindent(c) In the cases $\lieg{C}_{n}$ and $\lieg{D}_{n}$ the 
form $(...,...)$ induces a 
nondegenerate form on span$_{\mathbb{C}}\{e_{n-i},e_{n-i+1},..,e_{n}%
,e_{n+1},...,e_{n+i+1}\}$. 
 
\smallskip 
 
If $\mathfrak{p}$ is a standard parabolic subalgebra of one of 
these Lie algebras then it is given by a palindromic sequence of 
the form (A) 
$n_{1},...,n_{r-1},n_{r},n_{r},n_{r-1},...,n_{1}$ or (B) $n_{1},...,n_{r}%
,m,n_{r},...,n_{1}$ with $n_{i}>0$ and $2(n_{1}+...+n_{r})=q$ in 
the former case and in the latter in addition $m>0$ and 
$2(n_{1}+...+n_{r})+m=q$. We will use the corresponding grade on 
the Lie algebra. Each of these sequences delineates a block form 
for the $q\times q$ matrices. The elements of $\mathfrak{g}_{1}$ 
have block forms 
\[ 
x=\left[ 
\begin{array} 
[c]{cccccccc}%
0 & X_{1} & 0 & 0 & 0 & 0 & 0 & 0\\ 
\vdots & \vdots & \ddots & \vdots & \vdots & \vdots & \vdots & \vdots\\ 
0 & 0 & 0 & X_{r-1} & 0 & 0 & 0 & 0\\ 
0 & 0 & 0 & 0 & X_{r} & 0 & 0 & 0\\ 
0 & 0 & 0 & 0 & 0 & Y_{r-1} & 0 & 0\\ 
\vdots & \vdots & \vdots & \vdots & \vdots & \vdots & \ddots & \vdots\\ 
0 & 0 & 0 & 0 & 0 & 0 & 0 & Y_{1}\\ 
0 & 0 & 0 & 0 & 0 & 0 & 0 & 0 
\end{array} 
\right] 
\] 
for case (A) and%
\[ 
x=\left[ 
\begin{array} 
[c]{ccccccc}%
0 & X_{1} & 0 & 0 & 0 & 0 & 0\\ 
\vdots & \vdots & \ddots & \vdots & \vdots & \vdots & \vdots\\ 
0 & 0 & 0 & X_{r} & 0 & 0 & 0\\ 
0 & 0 & 0 & 0 & Y_{r} & 0 & 0\\ 
\vdots & \vdots & \vdots & \vdots & \vdots & \ddots & \vdots\\ 
0 & 0 & 0 & 0 & 0 & 0 & Y_{1}\\ 
0 & 0 & 0 & 0 & 0 & 0 & 0 
\end{array} 
\right] 
\] 
in case (B). Furthermore, for $i=1,..,r-1$ we look upon $X_{i}$ 
as a map of $V_{n_1+\dots+n_{i-1}+1,n_i}$ to 
$\check{V}_{n_1+\dots+n_{i-1}+1,n_i}$. 
and $Y_{i}=-X_{i}^{\#}$ with 
upper $\#$ indicating the adjoint map relative to the indicted 
perfect pairing. \ In case (A) $V_{n_1+...+n_{r-1}+1,n_{r}}=\check 
{V}_{n_{1}+...+n_{r-1}+1,n_{r}}$ so relative to the indicted 
non-degenerate forms indicated in (b) and (c) above we have 
$X_{r}^{\#}=-X_{r}$. \ In case (B) we 
look upon $X_{r}$ as a map from $V_{n_{1}+...+n_{r}+1,m}$ 
to $\check{V}_{n_{1}+...+n_{r}+1,m}$ and then $Y_{r}=-X_{r}^{\#}$. 
 
We will now begin a calculation of the Jordan canonical form of a 
generic element $x$ in $\mathfrak{g}_{1}$ as given above. \ As in 
the case A$_{n}$ we must calculate the rank of $x^{j}$ for 
$j=1,...,2r-1$ in case (A) and $j=1,...,2r$ in case (B). \ We need 
the following Lemma to handle the case (B) for C$_{n}$. \ Let 
$V_{1},...,V_{r},W$ be finite dimensional non-zero vector spaces 
and assume that we have a symplectic (i.e. non-degenerate skew 
symmetric bilinear form), $\left\langle ...,...\right\rangle $ on 
$W$. \ If $T:V_{i+1}\rightarrow V_{i}$ linear then define 
$T^{\ast}:V_{i}^{\ast }\rightarrow V_{i+1}^{\ast}$ (as usual) by 
$T^{\ast}\lambda=\lambda\circ T$. If $T:W\rightarrow V_{r}$ then 
we define $T^{\#}:V_{r}^{\ast}\rightarrow W$ by $\left\langle 
T^{\#}\lambda,w\right\rangle =\lambda(Tw)$. We leave the following 
as an exercise whose proof is similar to that of 
Lemma~\ref{lm:basic} and will be left to the reader. 
\begin{lm}\label{lm:basicC} 
Let $T_{i}:V_{i+1}\rightarrow V_{i}$ be linear for $i=1,...,r-1$ 
and let $T_{r}:W\rightarrow V_{r}$ be linear then the form on 
$V_{1}^{\ast}$ defined by 
\[ 
(\lambda,\mu)_{T}=\lambda(T\mu) 
\] 
with $T=T_{1}\cdots T_{r-1}T_{r}T_{r}^{\#}T_{r-1}^{\ast}\cdots 
T_{1}^{\ast} $ is skew symmetric. \ Let $X_{k}$ be the space of 
all skew symmetric forms of rank $k$ on $V_{1}^{\ast}$. Let $s$ be 
the minimum of $\{\dim V_{1},...,\dim V_{r},\dim W\}$ then if $s$ 
is even set $t=s$ if $s$ is odd set $t=s-1$ then 
the map%
 
\[ 
Hom(V_{2},V_{1})\times Hom(V_{3},V_{2})\times\cdots\times Hom(V_{r}%
,V_{r-1})\times Hom(W,V_{r})\rightarrow X_{k}%
\] 
given by $(T_{1},...,T_{r})\mapsto(...,...)_{T}$ with 
$T=T_{1}\cdots T_{r-1}T_{r}T_{r}^{\#}T_{r-1}^{\ast}\cdots 
T_{1}^{\ast}$ is onto. 
\end{lm} 
In case (A) the non-zero entries of $x^{j}$ are in the blocks in 
the $j$-th diagonal (i.e. has blocks $i,i+j$, so the diagonal is 
the $0$-th diagonal). 
 The entries in that diagonal are of three possible types: 
 
\smallskip 
 
\noindent1-A $X_{i}X_{i+1}\cdots X_{i+j-1}$ if $i\geq1$ and 
$j+i\leq r+1$ (the $i$-th entry). 
 
\smallskip 
 
\noindent2-A $X_{i}X_{i+1}\cdots X_{r}Y_{r-1}\cdots Y_{2r-i-j+1}$ 
if $i\geq1$, $2r>i+j>r+1$ and $i\leq r$ (the $i$-th entry). 
 
\smallskip 
 
\noindent3-A $Y_{r-i}Y_{r-i-1}\cdots Y_{r-i-j+1}$ if 
$i\geq1,i+j\leq r$ the ($r+i$-th entry). 
 
\smallskip 
 
In case (B) there are also three possibilities: 
 
\smallskip 
 
\noindent1-B $X_{i}X_{i+1}\cdots X_{i+j-1}$ if $i\geq1$ and 
$j+i\leq r+1$ (the $i$-th entry). 
 
\smallskip 
 
\noindent2-B $X_{i}X_{i+1}\cdots X_{r}Y_{r}\cdots Y_{2r-i-j+1}$ if 
$i\geq1$ and $2r>i+j>r+1$ and $i\leq r$ (the $i$-th entry). 
 
\smallskip 
 
\noindent3-B $Y_{r-i}Y_{r-i-1}\cdots Y_{r-i-j+1}$ if $i\geq0$ and 
$i+j\leq r$ the ($r+i+1$-th entry). 
 
\smallskip 
 
We will now compute the generic rank of each of the six types of 
matrices. 
\begin{prop}\label{prop:mainBCD} 
(i) In cases 1-A, and 1-B and $i+j\leq r$ the generic rank is 
$\min \{n_{i},...,n_{i+j}\}$. 
 
(ii) In case 1-A with $i+j=r+1$ the generic rank is $\min\{n_{i}%
,...,n_{r-1},s\}$ with $s=n_{r}$ if $n_{r}$ is even or the algebra 
is symplectic and $s=n_{r}-1$ if $n_{r}$ is odd and the algebra is 
orthogonal. 
 
(iii) In case 1-B and $i+j=r+1$ the generic rank is $\min\{n_{i}%
,...,n_{r},m\}$. (Here recall that $m$ is the block size of the 
middle block.) 
 
(iv) In case 2-A if $i=2r-i-j$ set 
$s=\min\{n_i,\dots, n_{2r-i-j}\}$. The generic rank is 
$s$ if $s$ is even and $s-1$ if $s$ is odd. 
If $i<2r-i-j$ set 
$s=\min\{n_{r},...,n_{2r-i-j+1}\}$ and 
$t=\min\{n_{i},n_{i+1},..,n_{2r-i-j}\}$, if $i>2r-i-j$ then set 
$s=\min \{n_{i},...,n_{r}\}$ and $t=\min\{n_{2r-i-j},...,n_{i}\}$ 
then there are two possibilities. 
 
\qquad a) If $s$ is even or if we are in the case $\lieg{C}_{n}$ 
then the generic rank is $\min\{s,t\}$. 
 
\qquad b) If $s$ is odd then the generic rank is $\min\{s-1,t\}$. 
 
(v) In case 3-A it is $\min\{n_{r-i+1},n_{r-i},...,n_{r-i-j+1}\}$. 
 
(vi) In case 3-B it is $\min\{n_{r-i+1},n_{r-i},...,n_{r-i-j+1}\}$ 
if $i>0$ and it is 
 
$\min\{m,n_{r},n_{r-i-1},...,n_{r-i-j+1}\}$ if $i=0$. 
 
(vii) In case 2-B it is $\min\{n_{i},...,n_{r},m,n_{r},n_{r-1},...,n_{2r-i-j}%
\}$ for orthogonal algebras. 
 
(viii) In case 2-B for symplectic algebras let $k$ be the minimum 
of $i$ and $2r-i-j$. \ Let $s=\min\{n_{k},...,n_{r},m\}$ if $s$ is 
even then set $t=s$ if 
$s$ is odd set $t=s-1$. \ Then the generic rank is $\min\{n_{i},...,n_{r}%
,m,n_{r},n_{r-1},...,n_{2r-i-j},t\}$. 
\end{prop} 
 
\begin{proof} 
The cases 1-A for $i+j\leq r$ ,3-A,1-B,3-B have already been 
observed in the study of $\lieg{A}_{n}$. \ In the case of 1-A with 
$i+j=r+1$ then in the symplectic case the generic rank of $X_{r}$ 
is $n_{r}$ in the orthogonal case it is $n_{r}$ if $n_{r}$ is even 
and $n_{r}-1$ if it is odd. We look at the case 2-B. This case for 
orthogonal algebras follows from Lemma 1 above, the similar fact 
that (in the notation of Lemma 1) the set 
\[ 
Z=\{XX^{T}|X\in Y_{p,q,r}\}=\{Z|Z^{T}=Z,Z\in 
M_{p}(\mathbb{C}),rk(Z)=r\} 
\] 
is irreducible and the last part of the proof of Lemma 1. \ For 
symplectic algebras we use the previous lemma and a similar 
argument. We are left with 2-A. In this case we must look at 
matrices of the form $XZX^{T}$ with the matrix $Z,$ $l\times l$ 
and $Z^{T}=-Z$ for the orthogonal cases and $MZ^{T}M=Z$ (with $M$ 
as above) in the symplectic case. In the orthogonal cases if $X$ 
is generic of rank $a$ and $Z$ is generic of rank $b$ (note $b$ is 
$l$ if $l$ is even and $l-1$ if $l$ is odd) then the rank of 
$XZX^{T}$ is $c=\min\{a,b\}$ if $c$ is even and it is $c-1$ if $c$ 
is odd. \ In the symplectic case the generic rank is the minimum 
of the generic ranks of $X$ and $Z$ (which must be $l$). 
\end{proof} 
We will use the result above to calculate the Jordan forms of 
specific generic elements in the next few sections using formula 
(*) above. 
 
\begin{cor}\label{cor:sameas} 
Let $x$ be a generic element in $\mathfrak{g}_{1}$ for the 
standard parabolic subalgebra $\mathfrak{p}$. If $\mathfrak{g}$ is 
of type $\lieg{C}_{n}$ and the corresponing sequence of block 
sizes is as in case (A) or if it is of type $\lieg{B}_{n}$ or 
$\lieg{D}_{n}$ then the Jordan form of $x$ is the same as that of 
$x^{\prime}$, a generic element in level 1 of the grade for the 
sequence of block sizes for $\mathfrak{sl}_{q}$ with $q=2n$ for 
types $\lieg{C}_n$ and $\lieg{D}_n$ and $2n+1$ for type 
$\lieg{B}_n$. 
\end{cor} 
 
We will conclude this section with some observations that will be 
used in the determination of the nice parabolic subalgebras that 
are given in terms of even nilpotent elements as in 
Lemma~\ref{lm:even}. 
 
Let $\liea{g}$ be a semisimple Lie algebra over $\CC$ and let 
$\liea{k}$ be a reductive subalgebra that is reductively embedded 
(that is the restriction of the adjoint representation of 
$\liea{g}$ to $\liea{k}$ is completely reducible).  Let $\liea{a}$ 
be a Cartan subalgebra of $\liea{k}$ contained in $\liea{h}$, a 
Cartan subalgebra of $\liea{g}$. Then we say that $\liea{k}$ is 
evenly embedded in $\liea{g}$ if the restriction to $\liea{a}$ of 
a root of $\liea{g}$ on $\liea{h}$ is in the root lattice of 
$\liea{k}$ in $\liea{a}^*$.  We note that this condition is 
independent of any choices made.  We note that if $\liea{k}$ is an 
embedded $\mathfrak{sl}_{2}$ with standard basis $\{e,f,h\}$ then 
$\liea{k}$ is evenly embedded if and only if $e$ is an even 
nilpotent element of $\liea{g}$. 
 
The following lemma is an obvious consequence of the definitions. 
 
\begin{lm} Let $\liea{k}$ be evenly embedded in the semisimple Lie 
algebra $\liea{g}$. Then $x \in \liea{k}$ is an even nilpotent 
element of $\liea{k}$ if and only if it is even as a nilpotent 
element of $\liea{g}$. 
\end{lm} 
 
We apply this result to the classical Lie algebras with the 
following easy Lemma that we will leave to the reader. 
\begin{lm} 
Standard embeddings of Lie algebras $\mathfrak{so}(n)$ and 
$\mathfrak{sp}(n)$ into $\mathfrak{sl}_{d}$ (with $d=n$ and $d=2n$ 
respectively) are even. 
\end{lm}

\section{The proofs of the main theorems.} 
In this section we will prove the theorems of section 1. For the 
exceptional types the results were proved using the computer 
program GAP.  Our proofs in the classical cases involve the 
computation of Jordan forms for generic elements as described in 
Proposition ~\ref{prop:jordA} and Proposition~\ref{prop:mainBCD}. 
\subsection{The case of $\lieg{A}_n$} 
In this subsection we will prove Theorem ~\ref{th:A}. We will 
first prove the necessity of the condition. We will use the 
following result 
 
\begin{thm}\label{th:centA} 
Let $n_{1}\geq n_{2}\geq...\geq n_{r}>0$ define the Jordan form of 
a nilpotent 
matrix in $M_{n}(\mathbb{C)}$ then the dimension of its centralizer is $%
{\textstyle\sum\limits_{i}} 
m_{i}^{2}$ with $m_{1}\geq m_{2}\geq...\geq m_{s}$ the dual 
partition. 
\end{thm} 
Suppose that there exist $i<j<k$ with $n_{i}>n_{j}<n_{k}$. \ We 
assert that 
this implies that there exist $l$ and $m$ such that $n_{l}<n_{l+1}%
=...=n_{i+m}>n_{l+m+1}$. \ Indeed, let $k-i$ be minimal subject to 
the condition that there exists $i<j<k$ such that 
$n_{i}>n_{j}<n_{k}$. If there exists $s$ with $i<s<j$ and 
$n_{s}\geq n_{i}$ then $n_{s}>n_{j}<n_{k}$ and $k-s<k-i$. \ Thus 
if $i<s<j$ then $n_{i}>n_{s}$. If $i<s<j$ and $n_{s}<n_{j}$ then 
$n_{i}>n_{s}<n_{j}$ and $j-i<k-i$. So we must have if $i<s<j$ then 
$n_{i}>n_{s}\geq n_{j}$. If for such an $s$ we have $n_{s}>n_{j}$ 
then we have $n_{s}>n_{j}<n_{k}$. \ Thus we have $n_{s}=n_{j}$. 
For $i<s<j$. We can use the same argument to complete the proof 
for the interval between $j$ and $k.$ 
 
\begin{lm} 
If $\mathfrak{g}$ is of type $\lieg{A}_n$ and if $\mathfrak{p}$ is 
the standard 
parabolic subalgebra of $\lieg{A}_n$, $n\geq1$ given by integers $n_{1}%
,...,n_{r}>0$ with $n_{1}+...+n_{r}=n+1$. Then if there exist 
$i<j<k\leq r$ with $n_{i}>n_{j}<n_{k}$ then $\mathfrak{p}$ is not 
nice. 
\end{lm} 
\begin{proof} 
The above remarks imply that we may assume that we have%
\[ 
n_{1}>n_{2}=...=n_{m}<n_{m+1}. 
\] 
Then the above proposition implies that if the Jordan canonical 
form of a generic element of $\mathfrak{g}_{1}$ has Jordan 
canonical form with $a_{j}$ blocks of size $j$ for $j=1,...,m+1$ 
then $a_{m+1}=n_{2}$ and $a_{i}=0$ for 
$1<i<m+1$ and $a_{1}=n_{1}+n_{m+1}-2n_{2}$. The corresponding partition is%
\[ 
m+1=...=m+1>1=1=...=1 
\] 
with $n_{2}$, $m+1$ and $n_{1}+n_{m+1}-2n_{2}$ one. The dual 
partition has one row of length $n_{1}+n_{m+1}-n_{2}$ and $m$ rows 
of length $n_{2}$. Thus the dimension of the centralizer of a 
generic element of $\mathfrak{g}_{1}$ in $M_{n+1}(\mathbb{C})$ is 
\begin{align*} 
(n_{1}+n_{m+1}-n_{2})^{2}+mn_{2}^{2}  &  =2(n_{1}-n_{2})(n_{m+1}-n_{2}%
)+n_{1}^{2}+n_{m+1}^{2}+(m-1)n_{2}^{2}\\ 
&  >n_{1}^{2}+n_{m+1}^{2}+(m-1)n_{2}^{2}=\dim\mathfrak{l}%
\end{align*} 
with $\mathfrak{l}$ the standard Levi factor of the corresponding 
parabolic subalgebra in $M_{n+1}(\mathbb{C})$. \ Thus 
Theorem~\ref{th:centA} implies that $\mathfrak{p}$ is not nice. 
\end{proof} 
In light of this lemma the condition of the main theorem for 
A$_{n}$ is necessary. \ We now prove that it is sufficient. We 
will prove the result by induction on $r$. \ If $r=1$ there is 
nothing to prove. The above lemma implies that we may assume that 
the $n_{i}$ satisfy one of the following 
systems of inequalities%
\[ 
n_{1}\leq...\leq n_{r}\underset{}{>0\qquad(i)}%
\]%
\[ 
n_{1}\geq...\geq n_{r}\underset{}{>0\qquad(ii)}%
\]%
\[ 
0<n_{1}\leq...\leq n_{r_{1}}>n_{r_{1}+1}\geq...\geq n_{r}\underset{}%
{>0\qquad(iii)}%
\] 
Using the Chevalley involution we see that case (i) is nice if and 
only if case (ii) is nice. \ But Proposition~\ref{prop:jordA} 
above implies that in case (ii) the partition corresponding to the 
Jordan canonical form of a generic element of $\mathfrak{g}_{1}$ 
is the dual partition to $n_{1}\geq...\geq n_{r}$. Thus the 
dimension of the centralizer of a generic element in 
$\mathfrak{g}_{1}$ in $M_{n+1}(\mathbb{C})$ is $%
{\textstyle\sum} 
n_{i}^{2}$. \ Thus the parabolic subalgebra is indeed nice in 
these cases. We now complete the proof by induction. \ If $r\leq2$ 
then we are in case (i) or (ii).\ \ We assume the result for 
$r-1\geq2$. \ We may assume that we are in case (iii). \ Applying 
the Chevalley automorphism we may assume that $n_{1}\geq n_{r}$. 
The formula for the Jordan form implies that if $i>r-1$ the 
contribution to the $a_{i}$ of $n_{r}$ cancel out. \ Let 
$a_{i}^{\prime}$ be the number of blocks of size $i$ for 
\[ 
n_{1}\leq...\leq n_{r_{1}}>n_{r_{1}+1}\geq...\geq n_{r-1}. 
\] 
Then $a_{i}^{\prime}=a_{i}$ for $i=2,...,r-2$ and $a_{r-1}^{\prime}%
=\min\{n_{1}.n_{r-1}\}$. We also note that $a_{r-1}=\min\{n_{1},n_{r-1}%
\}-n_{r}$ and $a_{r}=n_{r}$. This implies that 
$a_{1}^{\prime}=a_{1}$. Hence the dual partition for the $a_{i}$ 
with $i=1,...,r$ differs from the dual partition for the 
$a_{i}^{\prime}$ with $i=1,...,r-1$ by having one additional row 
which is of length $n_{r}$. \ The inductive hypothesis implies 
that the 
sum of the squares of the entries in the dual partition for the $a_{i}%
^{\prime}$ is $n_{1}^{2}+...+n_{r-1}^{2}$ and so the sum of the 
squares for the dual partition for the $a_{i}$ is 
$n_{1}^{2}+...+n_{r}^{2}$. \ This completes the proof. 
 
Actually we have proved that if we have a parabolic subalgebra in 
the case A$_{n}$ that corresponds to a unimodal sequence 
$n_{1},...,n_{r}$ and if $p$ is the partition corresponding to the 
Jordan form of a generic element of $\mathfrak{g}_{1}$ then the 
dual partition to $p$ is the $n_{i}$ put in decreasing order. 
 
The following result follows from the fact that a nilpotent 
element in $\liea{sl}_{n}$ is even if and only if all of the block 
lengths in its Jordan canonical form have the same parity. 
 
\begin{thm}\label{th:evenA} Let $\liea{p}$ be a nice standard parabolic 
subalgebra of $\mathfrak{sl}_{n}$ then it corresponds to an even 
nilpotent element if and only if its sequence of block lengths is 
palindromic. 
\end{thm} 
 
This result implies that the even nilpotent orbits on $\mathfrak{sl}_n$ 
are in one to one correspondence with the 
unimodal, palindromic, compositions of $n$ (recall that a 
composition of $n$ is a 
sequence $m_{1},...,m_{r}$ of positive integers such that $m_{1}+...+m_{r}%
=n$). \ Let $a_{n}$ denote the number of even nilpotent orbits 
($0$ is considered to be an even nilpotent orbit) in 
$\mathfrak{sl}_n$. It is not hard to show that this 
description implies 
\begin{lm} 
The generating function for the $a_{n}$, 
$\sum_{n\geq1}a_{n}q^{n}$, is 
\[\sum\limits_{j\geq1}\frac{q^{j}(1+q^{j})} {\prod\limits_{i=1}^{j} 
(1-q^{2i})}.\] 
\end{lm} 
 
We note that this combined with the usual classification of even 
nilpotent 
orbits yields%
\[\sum\limits_{j\geq1}\frac{q^{j}(1+q^{j})} {\prod\limits_{i=1}^{j} 
(1-q^{2i})}=\frac{1}{ \prod\limits_{i=1}^{\infty} 
(1-q^{2i+1})}+\frac{1}{ \prod\limits_{i=1}^{\infty} (1-q^{2i})} 
\] 
This sequence can be found in the Online Encyclopedia of Integer 
Sequences as number A096441. \ We note that Euler's identity (cf. 
[A] (2.2.5)) says that 
\[ 
1+ \sum\limits_{n=1} 
\frac{t^{n}}{\prod\limits_{i=1}^{n}(1-q^{i})}= \frac{1}{ 
\prod\limits_{i=0}^{\infty}(1-tq^{i})}. 
\] 
This implies the above equation also. However, the above method is 
purely combinatorial and implies the special cases $t=q$ and 
$q\rightarrow q^{2}$ and $t=q^{2}$ and $q\rightarrow q^{2}$.

\subsection{The case of $\lieg{C}_n$.} 
In this subsection we will prove Theorem \ref{th:C}. We will use 
the following result (cf.~\cite{cm}). 
\begin{thm}\label{thm:centralizer-dim-C} 
Let $\mathfrak{g}$ be of type $\lieg{C}$. Let $n_{1}\geq 
n_{2}\geq...\geq n_{r}>0$ be a partition corresponding the Jordan 
form of an element $x\in\mathfrak{g} $ then if $m_{1}\geq...\geq 
m_{s}$ is the dual partition then the dimension of 
$\mathfrak{g}^{x}$ is

\[\textstyle\sum\limits_{i} 
\frac{m_{i}^{2}}{2}+\frac{1}{2}|\{i|n_{i}\ odd\}|\] 
\end{thm}

Before we get down to the business at hand we will prove a useful 
result. 
 
\begin{lm}\label{lm:firstC} 
Let the standard parabolic subalgebra, $\mathfrak{p}$, for $C_{n}$ 
have block sizes $r$ of size $l$, $s$ of size $m$, and $r$ of size 
$l$ with $s$ and $l$ 
odd, $l<m$ and $n=rl+\frac{sm}{2}$ viz.%
\[ 
\left[ 
\begin{array} 
[c]{ccc}%
l\times l & \ast & \ast\\ 
0 & m\times m & \ast\\ 
0 & 0 & l\times l 
\end{array} 
\right]  . 
\] 
Then $\mathfrak{p}$ is nice if and only if $r=0,1$. 
\end{lm} 
\begin{proof} 
We will use the main result of the previous section. \ Let $x$ be 
a generic element in $\mathfrak{g}_{1}$. Let $r_{i}$ denote the 
rank of $x^{i}$. Our 
result implies that $r_{i}=2rl+(s-i)m$ for $i=1,...,s$, $r_{s+i}%
=(2r-2i)l+i(l-1)$ for $i=1,...,r$ and $r_{s+r+i}=(r-i)(l-1)$ for 
$i=1,...,r$. Thus using the method of section one to calculate 
Jordan forms we see that the Jordan form of a generic $x$ has 
$l-1$ blocks of size $2r+s$, $2$ blocks of size $r+s$ and $m-l-1$ 
blocks of size $s$. 
 
The dual partition has $s$ blocks of size $m$, $r$ blocks of size 
$l+1$ and $r$ blocks of size $l-1$. We therefore see that the 
dimension of the centralizer of $x$ is 
\[ 
\frac{sm^{2}}{2}+(l^{2}+1)r+\frac{m}{2}-1+\xi 
\] 
with $\xi=1$ if $r$ is even and $0$ if $r$ is odd. \ This number 
is $\frac{s-1}{2}m^{2}+\frac{m(m+1)}{2}+l^{2}r+r-1+\xi$. If $r$ is 
even this number is $r$ and if $r$ is odd it is $r-1$. The lemma 
follows. 
\end{proof} 
We will first show that a standard parabolic subalgebra is nice 
only if the corresponing sequence defining the A$_{2n-1}$ 
parabolic is unimodal. So suppose that we have a standard 
parabolic subalgebra $\mathfrak{p}$ corresponding to 
$n_{1},...,n_{r},n_{r},...,n_{1}$ in case (A) or 
$n_{1},...,n_{r},m,n_{r},...,n_{1}$ in case (B). We assume that 
the sequence is not unimodal and we derive a contradiction. \ We 
first consider case (A). Then Corollary ~\ref{cor:sameas} implies 
that the Jordan form in question is the same as the one for 
A$_{q}$ with $q=2(n_{1}+...+n_{r})$. \ Let $m_{1}\geq...\geq 
m_{s}$ be the corresponding dual partition. If the sequence 
$n_{1},...,n_{r},n_{r},...,n_{1}$ is not unimodal then we have%
\[ 
2(n_{1}^{2}+...+n_{r}^{2})<m_{1}^{2}+....+m_{s}^{2}%
\] 
by the main theorem for A$_{n}$. This implies that%
\[ 
n_{1}^{2}+...+n_{r}^{2}<\frac{(m_{1}^{2}+....+m_{s}^{2})}{2}\leq 
\dim\mathfrak{g}^{x}%
\] 
so the parabolic is not nice. 
 
We now look at case (B). We assume that the squence $n_{1},...,n_{r}%
,m,n_{r},...,n_{1}$ is not unimodal but the corresponding 
parabolic subalgebra 
is nice. Then there must be a subsequence of the form%
\[ 
n_{s}>n_{s+1},...,m,...,n_{s+1}<n_{s}%
\] 
or of the form%
\[ 
n_{r}>m<n_{r}. 
\] 
These correspond to subdiagrams of the diagram for the standard 
parabolic subalgebra. \ We assert that the corresponding parabolic 
subalgebra cannot be nice for the corresponding subalgebra of type 
C$_{n_{r}+\frac{m}{2}}$. Let $x$ 
be generic in $\mathfrak{g}_{1}$. Then according to Proposition~\ref{prop:rank-C} 
\[ 
rk(x)=2m, 
\]%
\[ 
rk(x^{2})=m, 
\]%
\[ 
rk(x^{3})=0. 
\] 
Thus the Jordan canonical form of $x$ has $2n_{r}-2m$ blocks of 
size $1$ and $m$ blocks of size $3$. \ Thus the dual partition is 
$2n_{r}-m\geq m=m$. So the dimension of $\mathfrak{g}^{x}$ is 
\begin{align*} 
&  \frac{(2n_{r}-m)}{2}^{2}+m^{2}+n_{r}+\frac{m}{2}\\ 
&  =(n_{r}-m)^{2}+n_{r}+n_{r}^{2}+\frac{m(m+1)}{2}. 
\end{align*} 
Since $n_{r}^{2}+\frac{m(m+1)}{2}$ is the dimension of the 
corresponding Levi factor we see that the corresponding parabolic 
is not nice. \ We can thus 
assume that we have%
\[ 
n_{s}>n_{s+1},...,m,...,n_{s+1}<n_{s}%
\] 
as a subsequence. \ We take $s$ to be maximal with respect to this 
condition. 
\ Then we must have%
\[ 
n_{s}>n_{s+1}\leq...\leq n_{r},m,n_{r}\geq...\geq n_{s+1}<n_{s}%
\] 
if one of the inequalities $n_{t}\leq n_{t+1}$ were strict then 
the type A result combined with the subdiagram theorem would imply 
that the parabolic is 
not nice. Thus we may assume that we have%
\[ 
n_{s}>n_{s+1}=..=n_{r},m,n_{r}=...=n_{s+1}<n_{s}%
\] 
the relationship between $n_{r}$ and $m$ cannot be $n_{r}>m$ since 
this yields a subdiagram that we have just seen cannot be nice. We 
are thus left with 2 
cases%
\[ 
a>b=...=b<a 
\] 
or%
\[ 
a>b=...=b<c>b=...=b<a 
\] 
where $a=n_{s}$ and in the former case $b=m$ in the latter case 
$n=n_{s+1}$ and $c=m$. We look at the latter case first. We assume 
first of all that $b$ is even. Our usual method of calculation 
(assuming that there are $r$ $b$'s 
and that $n=2rb+2a+m$) yields%
\[ 
\dim\mathfrak{g}^{x}-l=(a-b)(a+2c+1-3b)\text{.}%
\] 
If $b$ is odd then the previous lemma implies that $r=1$. \ Thus 
we are 
looking at%
\[ 
a>b<c>b<a. 
\] 
In this case one finds that the dual partition is $2a+c-2b\geq 
b+1=b+1>b-1=b-1$. We leave the algebra that shows that this case 
is not nice to the reader. 
 
At this point we have shown that unimodality is a necessity for 
niceness. The lemma above implies that it is not sufficient.  We 
will devote the rest of the section to the proof of Theorem 
~\ref{th:C}. Let $\mathfrak{p}$ be a standard parabolic of 
$\mathfrak{g}$ with a unimodal sequence of block lengths. Then if 
the number of blocks is even or all of the block lengths are even 
then Corollary ~\ref{cor:sameas} implies that the Jordan canonical 
form of a generic element in $\mathfrak{g}_{1}$ is the same as the 
Jordan form for the corresponding parabolic in $\lieg{A}_{2n-1}$. 
Since the sequence is unimodal and palindromic we see that the 
element corresponds to an even nilpotent element for $A_{2n-1}$ 
and hence for $\mathfrak{g}_{1}$. \ Thus in this case the 
parabolic is nice. \ We may thus assume that the number of 
blocks is odd and the sequence of block sizes is of the form%
\[ 
0<n_{1}\leq...\leq n_{r}\leq m\geq n_{r}\geq...\geq n_{1}>0.\overset{}%
{\quad\quad(1)}%
\] 
With at least one $n_{i}$ odd. Lemma ~\ref{lm:firstC} implies that 
the result is true for $r=1$. We will prove the theorem in this 
case by induction on $r$. We will prove that if we have a sequence 
as above that corresponds to a nice parabolic and if we extend it 
by adding blocks at either end of size 
$n_{0}$ yielding%
\[ 
0<n_{0}\leq n_{1}\leq...\leq n_{r}\leq m\geq n_{r}\geq...\geq 
n_{1}\geq n_{0}>0. 
\] 
Then the new parabolic subalgebra is nice if $n_{0}$ is even or if 
$n_{0}$ is odd and $n_{0}<n_{1}$. We will show that it is not nice 
if $n_{0}$ is odd and $n_{0}=n_{1}$. We will first prove the 
positive statements. 
 
Assume that $n_{0}$ is even. \ Let $x^{\prime}$ be a generic 
element of the first level of the grade for the parabolic 
subalgebra corresponding to the block sequence (1) and let $x$ be 
an element for the block sequence above. \ Let $R_{j}^{\prime}$ be 
the rank of $(x^{\prime})^{j}$ and let $R_{j}$ be the rank of 
$x^{j}$. \ Then since $n_{0}$ is even Proposition 
~\ref{prop:mainBCD} implies that 
\begin{align*} 
R_{j}  &  =R_{j}^{\prime}+2n_{0},1\leq j\leq2r+1,\\ 
R_{2r+2}  & 
=R_{2r+2}^{\prime}+n_{0},R_{2r+j}=R_{2r+j}^{\prime},j>2. 
\end{align*} 
This if $N_{j}$ and $N_{j}^{\prime}$ are the respective dimensions 
of the 
kernels of $x^{j}$ and $(x^{\prime})^{j}$ we have%
\begin{align*} 
N_{j}  &  =N_{j}^{\prime},1\leq j\leq2r+1,\\ 
N_{2r+2}  & 
=N_{2r+2}^{\prime}+n_{0},N_{2r+j}=N_{2r+j}^{\prime}+2n_{0},j>2. 
\end{align*} 
So if the multiplicities in the Jordan forms of $x$ and 
$x^{\prime}$ are 
respectively $M_{j}$ and $M_{j}^{\prime}$ then we have%
\begin{align*} 
M_{j}  &  =M_{j}^{\prime},1\leq j\leq2r,\\ 
M_{2r+1}  &  =M_{2r+1}^{\prime}-n_{0},M_{2r+2}=M_{2r+2}^{\prime}%
=0,M_{2r+3}=n_{0}. 
\end{align*} 
This implies that the number of odd blocks in the Jordan form is 
unchanged and that the only change in the dual partition is the 
addition of $2$ blocks of size $n_{0}$. Thus as in the case of 
type A the parabolic subalgebra is nice in this case. 
 
We now consider the case when $n_{0}$ is odd and $n_{0}<n_{1}$. \ 
The argument is similar. \ We will use the same notation. \ We 
have 
\begin{align*} 
R_{j}  &  =R_{j}^{\prime}+2n_{0},1\leq j\leq2r+1,\\ 
R_{2r+2}  & 
=R_{2r+2}^{\prime}+n_{0}-1,R_{2r+j}=R_{2r+j}^{\prime},j>2. 
\end{align*} 
Doing the obvious subtraction we have%
\begin{align*} 
N_{j}  &  =N_{j}^{\prime},1\leq j\leq2r+1,\\ 
N_{2r+2}  & 
=N_{2r+2}^{\prime}+n_{0}+1,N_{2r+j}=N_{2r+j}^{\prime}+2n_{0},j>2. 
\end{align*} 
This yields 
\begin{align*} 
M_{j}  &  =M_{j}^{\prime},1\leq j\leq2r,\\ 
M_{2r+1}  &  =M_{2r+1}^{\prime}-n_{0}-1,M_{2r+2}=M_{2r+2}^{\prime}%
+2,M_{2r+3}=n_{0}-1. 
\end{align*} 
We therefore see that we have reuced the number of odd blocks by 
$2$ and added $2$ blocks to the dual partition of siizes $n_{0}+1$ 
and $n_{0}-1$. This completes the proof in this case. 
 
Finally we assume that $n_{0}$ is odd and $n_{0}=n_{1}$. \ We may 
assume that $n_{1}<n_{2}$ by the inductive hypothesis. We will 
show that in this case the parabolic subalgebra is not nice and 
the proof will be complete. \ This time 
we have%
\begin{align*} 
R_{j}  &  =R_{j}^{\prime}+2n_{0},1\leq 
j\leq2r,R_{2r+1}=R_{2r+1}^{\prime 
}+2(n_{0}-1),\\ 
R_{2r+2}  & 
=R_{2r+2}^{\prime}+n_{0}-1,R_{2r+j}=R_{2r+j}^{\prime},j>2. 
\end{align*} 
So%
\begin{align*} 
N_{j}  &  =N_{j}^{\prime},1\leq j\leq2r,N_{2r+1}=N_{2r+1}^{\prime}+2,\\ 
N_{2r+2}  & 
=N_{2r+2}^{\prime}+n_{0}+1,N_{2r+j}=N_{2r+j}^{\prime}+2n_{0},j>2. 
\end{align*} 
Hence%
\begin{align*} 
M_{j}  &  =M_{j}^{\prime},1\leq j\leq2r-2,M_{2r-1}=M_{2r-1}^{\prime}%
,M_{2r}=M_{2r}^{\prime}-2=0\\ 
M_{2r+1}  & 
=M_{2r+1}^{\prime}-n_{0}+1=0,M_{2r+2}=M_{2r+2}^{\prime 
}=0,M_{2r+3}=n_{0}-1. 
\end{align*} 
Here we see that there is a net change of two more odd blocks in 
the Jordan form. \ In the dual partition we remove one row of 
size, $n_{0}-1$ and replace it with one of size $n_{0}+1$, we then 
add $2$ rows of length $n_{0}-1$ to complete it. \ Thus the 
dimension of the centralizer for $x$ is that of the one 
for $x^{\prime}$ plus%
\[ 
\frac{(n_{0}+1)^{2}}{2}-\frac{(n_{0}-1)^{2}}{2}+(n_{0}-1)^{2}+1=n_{0}^{2}+2. 
\] 
The inductive hypothesis now implies that the dimension of the 
centralizer of $x$ is the dimension of a Levi factor of 
$\mathfrak{p}$ plus $2$ and the result follows. 
 
Using Theorem \ref{th:evenA} we note that the methods of this 
section allow us to describe the standard  nice parabolic 
subalgebras coming from even nilpotent elements in the case of 
$\lieg{C}_{n}$ and thereby give a parametrization of the even 
nilpotent orbits. 
 
\begin{lm} 
The even nilpotent orbits of a Lie algebra of type $\lieg{C}_{n}$ 
are in one two one correspondence with the unimodal, palindromic 
compositions of $2n$ which have the property that if the number of 
parts is odd then all the parts are even. 
\end{lm} 
 
\subsection{The case of $\lieg{B}_n$ and $\lieg{D}_n$.} 
 
In this subsection we will describe the standard nice parabolic 
subalgebras in the case when $\mathfrak{g}$ is an orthogonal Lie 
algebra that is of type B or D and thereby prove Theorems 
\ref{th:B} and \ref{th:D}. We will use the following result 
(cf.~\cite{cm}). 
 
\begin{thm} 
Let $\mathfrak{g}$ be of type $\lieg{B}$ or $\lieg{D}$. Let 
$n_{1}\geq n_{2}\geq...\geq n_{r}>0$ be a partition corresponding 
the Jordan form of an element $x\in\mathfrak{g}$ then if 
$m_{1}\geq...\geq m_{s}$ is the dual partition then the dimension 
of $\mathfrak{g}^{x}$ is 
\[\textstyle\sum\limits_{i} 
\frac{m_{i}^{2}}{2}-\frac{1}{2}|\{i|n_{i}\ odd \}|.\] 
\end{thm} 
We will prove Theorems ~\ref{th:B} and ~\ref{th:D} simultaneously 
in the following form. 
 
\begin{thm} 
Let $\mathfrak{p}$ be a standard parabolic subgroup of 
$\mathfrak{g}$ we consider 4 conditions on the sequence of its of 
block lengths. 
 
1) It is a unimodal sequence with an odd number of elements. 
 
2) It is a unimodal sequence of even length and each odd block 
length occurs at most twice. 
 
3) The sequence is of the form 
\[ 
n_{1}\leq....\leq n_{r}>n_{r}-1=...=n_{r}-1<n_{r}\geq...\geq n_{1}%
\] 
with the total number of blocks odd and $r=1$ or $n_{r-1}<n_{r}$. 
 
4) The sequence is of the form above with the total number of 
blocks even, the odd ones occuring at most twice and $n_{r}$ must 
be odd and at least $3$. 
 
Then $\mathfrak{p}$ is nice if and only if the sequence satisfies 
one of the above conditions. 
\end{thm} 
\medskip The proof of this result will take up the rest of the 
section. We will first prove that the shapes of the sequences are 
necessary conditions for niceness. \ So assume that we a nice 
$\mathfrak{p}$ with block lengths $m_{1},...,m_{s}$ that is not 
unimodal. Since $m_{j}=m_{s+1-j}$ we see that\ If $m_{i}\leq 
m_{i+1}$ for all $i\leq\frac{s}{2}-1$ if $s$ is odd and 
$i\leq\frac{s-1}{2}$ if $s$ is odd then the sequence is unimodal. 
\ Thus we may assume that there is $j$ such that $m_{j}>m_{j+1}$ 
and $j$ is in the indicated range. If we have $m_{j}>m_{k}<m_{l}$ 
with $j<k<l<u$ and $u=\frac{s}{2}$ if $s$ is even, 
$u=\frac{s-1}{2}$ if $s$ is odd then we would have a non-nice type 
type A subdiagram of the diagram of $\mathfrak{p}$. Taking the $j$ 
with $m_{j}>m_{j+1}$ in the indicted range to be as big as 
possible we see that there must be central subsequences of one of 
the following forms (even and odd correspond to the number of 
blocks, thus (1) in the even cases means that there are an even 
number of $b$'s): 
 
Even Cases:%
 
\[ 
(1)\qquad\ a>b=...=b<a. 
\]%
\[ 
(2)\qquad\ a>b=...=b<c=c>b=...=b<a. 
\]

Odd Cases:%
\[ 
(1)\qquad\ a>b=...=b<a. 
\]%
\[ 
(2)\qquad\ a>b=...=b<c>b=...=b<a. 
\] 
We will show that the diagrams of type (2) are never nice and 
those of type (1) are nice if and only if and $a=b+1$ and if there 
are an even number of blocks then $b$ is even. We first look at 
the case of $b$ even we are in the even type (1) case so there are 
$2l$, $b$'s. Then applying Proposition~\ref{prop:rank-BD} we 
see that if $x$ is a generic element in the degree 
$1$ part of the grade and if $R_{i}$ is the rank of $x^{i}$ then%
 
\[ 
R_{1}=(2l+1)b,R_{2}=2lb,R_{3}=(2l-1)b,...,R_{2l+1}=b,R_{j}=0,j>2l+1. 
\] 
The dimension of $\ker x^{i}$ will be denoted $N_{i}$. We find 
that 
\[ 
N_{1}=2a-b,N_{2}=2a,N_{3}=2a+b,...,N_{2l+1}=2a+(2l-1)b,N_{j}=2a+2lb,j>2l+1. 
\] 
From this we find that the Jordan canonical form has $2(a-b)$ 
blocks of size $1$ and $b$ blocks of size $2l+2$  This implies 
that the Jordan form of $x$ is $2a-b,b,...,b$ with $2l+1$ $b$'s. 
The sequence of block lengths has 
$2(a-b)$ odd parts thus%
\[ 
\dim\mathfrak{g}^{x}=\frac{(2a-b)^{2}}{2}+\frac{2l+1}{2}b^{2}-a+b. 
\] 
If $\mathfrak{m}$ is the standard Levi factor of $\mathfrak{p}$ 
then 
$\dim\mathfrak{m}=a^{2}+lb^{2}$. A direct calculation yields that%
\[ 
\dim\mathfrak{g}^{x}-\dim\mathfrak{m}=(a-b)(a-b-1).%
\] 
This says that in this case the parabolic subalgebra is nice if 
and only if $a=b+1$. \ We next look at type 1 with an odd number 
of blocks then the argument is identical (since in the odd case 
the parity of $b$ plays no role. \ We conclude the analysis of the 
type one cases by looking at the even case with with $b$ odd and 
showing that this cannot correspond to a nice parabolic 
subalgebra. To simplify the calculation we first show that the 
case when we have blocks $b,b,b,b$ with $b$ odd is not nice. 
Indeed, in this case we have 
\[R_1 = 3b-1,R_2 = 2b-2,R_3 = b-1,R_4 = 0.\] 
Thus \[M_1 = 0,M_2 = 2, M_3 = 0, M_4 = b-1.\] If $b=1$ then the 
dual partition is $2,2$ so $\dim \liea{g}^x - \dim \liea{m} = 2$. 
If $b > 1$ then the dual partition is $b+1,b+1,b-1,b-1$ and we 
have $\dim \liea{g}^x - \dim \liea{m} = 2$. This proves our 
assertion. If we use Lemma ~\ref{lm:subdiag} then we are reduced 
to proving that the parabolic with block lengths $a > b=b < a$ is 
not nice if $b$ is odd. In this case we have 
\[R_1 = 3b-1,R_2 = 2b-2,R_3 = b-1,R_4=0.\] 
This yields \[M_1 = 2a-2b,M_2 = 2,M_3 = 0,M_4 = b-1.\] If $b=1$ 
then the dual partition is \[2a-b+1,b+1,b-1,b-1.\] Hence we have 
\[\dim \liea{g}^x - \dim \liea{m} = (a-b)^2+a-b+2 > 2.\] This 
proves the result in this case. 
 
\ We are thus left with analyzing the type (2) cases. When there 
are an even number of parts then there are $4$ cases (i) $b,c$ 
even,(ii) $c$ even, $b$ odd,(iii) $c$ odd,$b$ even and (iv) $b,c$ 
odd. We will not bore the reader by going through all of the $4$ 
type (2) even cases and the one type (2) odd case. \ We will do 
enough that the reader should have no trouble checking the missing 
cases. Recall that we are trying to prove that all of these cases 
are not nice. \ We look at the odd case (2) first. In this case we 
have 
\[ 
R_{j}=(2l+3-j)b. 
\] 
>From this we see that the corresponding Jordan form has $2a-3b+c$ 
blocks of size $1$ and $b$ blocks of size $2l+3$. The dual 
partition is one block of 
size $2a-2b+c$ and $2l+2$ blocks of size $b$. Thus%
\[ 
\dim\mathfrak{g}^{x}=\frac{(2a-2b+c)^{2}}{2}+(l+1)b^{2}-a+b-\frac{c}{2}. 
\] 
The dimension of the corresponding Levi factor is 
$a^{2}+lb^{2}+\frac {c(c-1)}{2}$. One checks easily that the 
dimension of $\mathfrak{g}^{x}$ is strictly larger than the 
dimension of the Levi factor. 
 
We now look at the even cases. We first look at the case when $b$ 
and $c$ are 
even and the number of $b$'s is $2l$ (thus the total length is $2l+4$). Then%
\[ 
R_{1}=(2l+2)b+c,\ R_{i}=(2l+4-i)b,\ 2\leq i\leq2l+4. 
\] 
This leads to a Jordan form with $2a-2b$ ones,$c-b$ twos and $b$ 
blocks of size $2l+4$. The dual partition has one block of size 
$2a-2b+c$, one block of size $c$ and $2l+2$ of size $b$. \ The 
dimension of a Levi factor is 
$a^{2}+lb^{2}+c^{2}$ and the dimension of $\mathfrak{g}^{x}$ is%
\[ 
\frac{(2a-2b+c)}{2}+\frac{(c-b)}{2}+(l+1)b^{2}-a+b 
\] 
the difference is $(a-b)(a+2c-3b-1)\geq2(a-b)>0$. Thus the 
parabolic subalgebra is not nice as asserted. If $c$ is odd and 
$b$ is even then the only change in the ranks is 
$R_{1}=(2l+2)b+c-1$. We leave it to the reader to show that the 
parabolic subalgebra is not nice in this case. \ We will now 
look at $b,c$ odd. Then we have%
\begin{align*} 
R_{1}  & 
=(2l+2)b+c-1,\quad R_{2}=(2l+2)b, \\ 
R_{j} & =(2l-2(j-3))b+(j-2)(b-1),\ 2<j\leq 
l+3,\\ 
R_{j}  &  =(2l+4-j)(b-1),\ l+4\leq j\leq2l+4. 
\end{align*} 
From this we find that the corresponding Jordan form has $2a-2b+2$ 
blocks of size $1$, $c-b-2$ blocks of size $2$, $2$ blocks of size 
$l+3$ and $b-1$ blocks of size $2l+4$. \ The number of odd blocks 
is $2a-2b+2+2\xi$ with $\xi=1$ if $l$ is even and $\xi=0$ if $l$ 
is odd. \ The dual partition has one block of size $2a+c-2b+1$, 
one block of size $c-1$, $l+1$ blocks of size $b+1$ 
and $l+1$ blocks of size $b-1$. We therefore have%
\begin{align*} 
\dim\mathfrak{g}^{x}  &  =\frac{(2a-2b+c+1)^{2}}{2}+\frac{(c-1)^{2}}%
{2}+(l+1)(b^{2}+1)-a+b-\xi\\ 
&  =2a^{2}-4ab+2b^{2}+2(a-b)(c+1)+c^{2}+1+lb^{2}+b^{2}+l-a+b-\xi\\ 
&  =a^{2}-4ab+3b^{2}+2(a-b)c+(1-\xi)+l+a^{2}+lb^{2}+c^{2}\\ 
&  =(a-b)^{2}+2(a-b)(c-b)+(1-\xi)+l+a^{2}+2lb+c^{2}. 
\end{align*} 
Since the dimension of the corresponding Levi factor is 
$a^{2}+2lb+c^{2}$ the corresponding parabolic subalgebra is not 
nice. \ If $c$ is even and $b$ is odd the only change in the ranks 
is $R_{1}=(2l+2)b+c$. We leave this case for the reader. 
 
We conclude that if a parabolic subalgebra is nice and the 
corresponding sequence of block sizes not unimodal then it must 
contain the following in the middle 
\[ 
b+1>b=..=b<b+1 
\] 
with $b$ even if the number of terms is even. \ We now prove that 
under these 
conditions on $b$ and the number of blocks the parabolic corresponding to%
\[ 
a\geq b+1>b=..=b<b+1\leq a 
\] 
is not nice. \ Let $s$ be the number of $b$'s, One checks that 
\[ 
R_{1}=2(b+1)+(s+1)b,R_{j}=(s+4-j)a,2\leq j\leq s+4. 
\] 
This implies that the corresponding Jordan form one block of size 
$2(a-b-1)$ (unless this number is $0$) has $2$ blocks of size $2$ 
and $b$ blocks of size $s+4$. Thus the dual partition is one block 
of size $2a-b$ one of size $b+2$ (notice that these are the same 
if $a=b+1)$ and $b$ with multiplicity $s+2$. 
Thus the dimension of $\mathfrak{g}^{x}$ is $\frac{(2a-b)^{2}}{2}%
+\frac{(b+2)^{2}}{2}+\frac{s+2}{2}b^{2}-(a-b-1)-\frac{b}{2}\xi$ 
with $\xi=0$ if $s$ is even and $\xi=1$ if $s$ is odd. The 
dimension of the corresponding Levi factor is 
$2(b+1)^{2}+\frac{s}{2}b^{2}-\frac{b}{2}\xi$. The difference is 
$2+(a-b)(a-b-1)$ so the parabolic subalgebra is not nice as 
asserted. \ We leave it to the reader to check that if $a$ is as 
above and $a<b+1$ then the parabolic subalgebra is nice. 
 
We will prove the theorem by induction on the number of elements 
in the sequence of block lengths. We will first consider the 
unimodal case we will show that if the number of blocks is odd 
then the corresponding parabolic subalgebra is nice. If the number 
of blocks is even then the corresponding parabolic subalgebra is 
nice if and only if each odd block occurs at most twice. \ We will 
prove this by induction on the number of blocks. If the number of 
blocks is one then the corresponding parabolic subalgebra is the 
whole Lie algebra and the result is obvious. If the number of 
blocks is $2$ then the corresponding parabolic subalgebra has 
commutative nilradical so is nice. Assume that we have proved the 
result for $r\geq1$ blocks. If we take 
such a unimodal nice parabolic subalgebra with sequence%
\[ 
0<n_{1}\leq...\geq n_{1}>0 
\] 
and we adjoin a block at both ends a block of size $n_{0}\leq 
n_{1}$ then we will show that if $r$ is odd then the corresponding 
parabolic subalgebra is nice and if $r$ is even then it is nice if 
$n_{0}$ is even or if it is odd it is nice if and only if 
$n_{0}<n_{1}$. \ We first assume that $r$ is odd or $n_{0}$ is 
even. We will use arguments similar to those for the case of type 
C. We will use the primed notation in the same way as in the 
previous section. \ We find that 
\begin{align*} 
R_{i}  &  =R_{i}^{\prime}+2n_{0},\ 1\leq i\leq r,\\ 
R_{r+1}  &  =R_{r+1}^{^{\prime}}+n_{0},\ R_{r+2}=R_{r+2}^{\prime}. 
\end{align*} 
Thus the multiplicities in the Jordan form are 
\[ 
M_{i}=M_{i}^{\prime}\text{ for }i\leq r-1, 
\]%
\[ 
M_{r}=M_{r}^{\prime}-n_{0},\ M_{r+1}=M_{r+1}^{\prime}=0,\ M_{r+2}=n_{0}. 
\] 
This implies that the corresponding dual partition is obtained by 
adjoining two rows of length $n_{0}$ to the original partition. \ 
Thus we have added $n_{0}^{2}$ to the dimension of the centralizer 
and the same to the dimension of the Levi factor so the parabolic 
subalgebra is nice in this case. We now assume that 
$n_{0}$ is odd and $n_{0}<n_{1}$. \ In this case we have%
\begin{align*} 
R_{i}  &  =R_{i}^{\prime}+2n_{0},\ 1\leq i\leq r,\\ 
R_{r+1}  & 
=R_{r+1}^{^{\prime}}+n_{0}-1,\ R_{r+2}=R_{r+2}^{\prime}=0. 
\end{align*} 
The multiplicities for the Jordan form are given by%
\[ 
M_{i}=M_{i}^{\prime}\text{ for }i\leq r-1, 
\]%
\[ 
M_{r}=M_{r}^{\prime}-n_{0}-1,\ M_{r+1}=M_{r+1}^{\prime}+2=2,\ M_{r+2}=n_{0}-1. 
\] 
This time we change the dual partition by adjoining a row of 
length $n_{0}-1$ and one of length $n_{0}+1$. In the Jordan form 
we have added to rows of odd length. So we find that the parabolic 
subalgebra is nice since the dimension of the Levi factor is the 
same as the centralizer. 
 
To complete the unimodal case we must show that if $r$ is even 
$n_{0}$ is odd and $n_{0}=n_{1}$ then the parabolic subalgebra is 
not nice. We note that by 
the inductive hypothesis we must have $n_{1}<n_{2}$. As above we have%
\begin{align*} 
R_{i}  &  =R_{i}^{\prime}+2n_{0},\ 1\leq i\leq r-1,\ R_{r}=R_{r}^{\prime}%
+2(n_{0}-1),\\ 
R_{r+1}  & 
=R_{r+1}^{^{\prime}}+n_{0}-1,\ R_{r+2}=R_{r+2}^{\prime}=0. 
\end{align*} 
This yields multiplicities in the Jordan form:%
\begin{align*} 
M_{i}  &  =M_{i}^{\prime}\text{ for }i\leq r-2,\ M_{r-1}=M_{r-1}^{\prime}-2=0,\\ 
M_{r}  &  =M_{r}^{\prime}-n_{0}+3=2,\ M_{r+1}=0,\ M_{r+2}=n_{0}-1. 
\end{align*} 
The upshot is two rows of length $r-1$ are lost and the original 
last row in the dual partition is replaced by a row of length 
$n_{0}+1$ and one row is added of length $n_{0}-1$. \ The upshot 
is that the dimension of the centralizer is two more than that of 
the Levi factor. We have completed the argument in this case. 
 
We now look at the non-unimodal case. Here we have seen that the 
sequence must 
contain%
\[ 
b+1>b=...=b<b+1 
\] 
with $b$ even if the number of terms is even or%
\[ 
a<b+1>b=...=b<b+1>a 
\] 
with $b$ even if the number of terms is even. Now using the 
results for type A we see that if we have a general nonunimodal 
sequence and not just of the 
above two forms it must be of the form%
\[ 
n_{1}\leq...\leq n_{r}\leq a<b+1>b=...=b<b+1>a\geq n_{r}\geq...\geq n_{1}%
\] 
with $b$ even if the number of blocks is even. \ Now exactly the 
same argument as in the unimodal case proves that the conditions 
3),4) of the theorem are necessary and sufficient. 
 
Using Theorem \ref{th:evenA} we note that the methods of this 
section allow us to describe the standard nice parabolic 
subalgebras coming from even nilpotent elements in the case of 
$\mathfrak{so}(n)$ and thereby give a parametrization of the even 
nilpotent orbits. 
 
\begin{lm} 
The even nilpotent orbits of $\mathfrak{so}(n)$ are in one-to-one 
correspondence with the unimodal, palindromic compositions of $n$ 
which have the property that if the number of parts is even then 
all the block sizes are even. 
\end{lm}

\subsection{The exceptional cases.} In this subsection we 
will explain how to derive the tables for the exceptional Lie 
algebras.  In the case of $\lieg{G}_2$, Lemma~\ref{lm:even} and 
Corollary~\ref{cor:dim} already give the classification of nice 
parabolic subalgebras. For $\lieg{F}_4$ there are eight parabolic 
subalgebras that are given by a TDS and seven parabolic 
subalgebras that are not nice by the dimension criterion 
Corollary~\ref{cor:dim}. The only remaining case is $(1,0,0,1)$. 
One uses Corollary~\ref{cor:dist} to show that this parabolic 
subalgebra is not nice: The corresponding Levi factor has 
dimension $12$. The only nilpotent element in the Bala-Carter 
table (cf.~\cite{ca}, p. 401) that has a $12$-dimensional 
centralizer satisfies $\dim\liea{m}=\dim\liea{g}_1$, i.e. is a 
distinguished nilpotent. But the grading associated to $(1,0,0,1)$ 
gives a much smaller dimension for $\liea{g}_1$. Hence, this 
parabolic subalgebra is not nice. 
 
\begin{re} 
Observe that for $\lieg{G}_2$ and $\lieg{F}_4$, the nice parabolic 
subalgebras are exactly those given by an even TDS. 
\end{re} 
 
The most complicated case is $\lieg{E}_n$. We describe the steps 
that led to the classification (cf. list in 
Section~\ref{sec:one}). 
 
First step: Use the dimension criterion (Corollary~\ref{cor:dim}). 
There remain $37$, $46$ resp. $40$ parabolic subalgebras of 
$\lieg{E}_6$, $\lieg{E}_7$, resp. of $\lieg{E}_8$. We used 
Mathematica to check the dimensions of the graded parts. 
 
In a second step, one identifies all parabolic subalgebras that 
are given by a TDS ($10$ for $\lieg{E}_6$, $24$ for $\lieg{E}_7$ 
and $46$ for $\lieg{E}_8$). They can be found using the tables of 
nilpotent orbits in~\cite{ca}, pp. 402-407. 
 
Thirdly, the subdiagram result (Lemma~\ref{lm:subdiag}) gives 
sixteen more nice parabolic subalgebras for $\lieg{E}_6$. It also 
excludes a couple of parabolics for all $\lieg{E}$-types that have 
a bad subdiagram (of type $\lieg{D}_n$, $n=4,5,6,7$). 
 
After that, there remain $30$ parabolic subalgebras of 
$\lieg{E}$-types. For all these parabolic subalgebras, one 
calculates the dimension of the Levi factor. By 
Theorem~\ref{thm:centr}, this is the same as the dimension of the 
centralizer for a Richardson element. In a case-by-case study, one 
checks whether there is a nilpotent element in $\liea{g}_1$ whose 
stabilizer has the required dimension. This can be done using the 
Computer Algebra Program GAP: 
 
\begin{re} 
The GAP-program uses the function `random' to generate a generic 
element $x\in\liea{g}_1$, namely it produces an element 
$x=\sum_{\alpha\in\liea{g}_1}c_{\alpha}X_{\alpha}$ where all 
$c_{\alpha}$ are nonzero. Generic means that the dimension $\dim 
\liea{m}^x$ of the centralizer of $x$ in the Levi factor 
$\liea{m}$ equals $\dim\liea{m}-\dim\liea{g}_1$. 
 
The program tests whether $x$ is generic and then calculates the 
dimension $\dim\liea{g}^x$ of the centralizer of $x$ in 
$\liea{g}$. By Theorem~\ref{thm:centr}, $\liea{p}$ is nice if and 
only if this dimension equals $\liea{m}$. 
\end{re} 
\begin{re} 
We came to the same list by hand. This involved studying the roots 
of $\liea{g}_1$ using a variety of different methods. 
\end{re} 
The complete list of nice parabolic subalgebras in type $\lieg{E}$ 
consists of $30$ parabolic subalgebras of $\lieg{E}_6$, $29$ of 
$\lieg{E}_7$ and $29$ of $\lieg{E}_8$, most of those come from 
an even TDS. The next table lists all nice parabolic subalgebras 
that are not given by an even TDS nor by a subdiagram of a 
parabolic given by an even TDS. We use the Bourbaki numbering 
of simple roots. 
\[ 
\begin{array}{ccc} 
(1,1,0,0,1,0) & (1, 1, 0, 0, 1, 0, 1) & (0,0,1,0,0,0,1,0)\\ 
(1,1,0,0,0,0) & (1, 1, 0, 0, 0, 0, 1) & \\ 
(0,1,1,0,0,1) & (0, 1, 1, 0, 0, 1, 1) & \\ 
(0,1,0,0,0,1) & (0, 0, 1, 0, 0, 0, 1) & \\ 
 & (0, 0, 0, 0, 1, 0, 1) & 
\end{array} 
\]

%
\section{Nilpotent elements, a normal form}\label{sse:X0} 
%
%
In this section, we introduce a normal form for Richardson 
elements in $\liea{g}_1$. 
We will not give any proofs given that the paper is already 
long. 
In a separate paper (\cite{ba}), 
the first named author is discussing 
this normal form for the nilpotent elements in more detail 
and explains the weight structure of the representation 
$\liea{g}_1$ of $\liea{m}$. 
 
Let $\liea{p}\subset\liea{g}$ be a nice 
parabolic subalgebra of a classical Lie algebra. We will show how 
to choose a Richardson element in $\liea{g}_1$ that projects into 
only the root spaces for the roots in a small set $S_1$. For 
each type of a parabolic subalgebra $\liea{p}$ in type 
$\lieg{A,B,C,D}$ we give a recipe that picks this set of 
$\liea{g}_1$.

We view the classical Lie algebras in the corresponding 
$\liea{gl}_N$ ($N=n+1,2n+1,2n,2n$ for types $\lieg{A,B,C,D}$). If 
$\alpha_{ij}$, $i<j$, is a positive root for $\liea{gl}_N$, then 
the standard matrix with a $1$ in the $i,j$ position and zeros 
everywhere else, $E_{ij}$, is a basis for the root space 
$\liea{g}_{\alpha_{ij}}$. 
 
Recall that the sequence of the block lengths 
of the standard Levi factor are of the form 
$(a_1,\dots,a_{r+1})$ for $\liea{sl}_{n+1}$. For 
the other classical Lie algebras, they are 
of the form $(a_1,\dots,a_r,a_r,\dots,a_1)$ 
(even number of blocks, (A)-case) or of the form 
$(a_1,\dots,a_r,a_{r+1},a_r,\dots,a_1)$ 
(odd number of blocks, (B)-case) respectively. 
These sequences give rise to a block decomposition 
of the first super-diagonal (cf. discussion above 
Lemma~\ref{lm:basicC}). We thus obtain $r$ subsets 
of the positive roots. More precisely,  the roots of 
$\liea{g}_1$ are divided into $r$ subsets. 
 
We can think of the roots of such a subset as filling 
out a 
rectangle of size $a_i\times a_{i+1}$, $1\le i\le r-1$ 
(resp. $i\le r$ in the (B)-case) 
and of  size $a_r\times a_r$ in the (A)-case. 
We denote these rectangles by $R_{i,i+1}$ ($i\le r$). 
Note that the central one in the (A)-case is a square matrix 
of size $a_r\times a_r$ that is symmetric if the Lie algebra 
is symplectic and skew-symmetric if $\liea{g}$ is an orthogonal 
Lie algebra. 
The entries of a rectangle $R_{i,i+1}$ will be labeled by 
$(k,l)$, $1\le k\le a_{i+1}$, $1\le l\le a_i$, starting 
with $(1,1)$ in the lower left corner and ending with 
$(a_{i+1},a_i)$ in the upper right corner. 
Every such entry $(k,l)$ corresponds uniquely to 
a positive root $\alpha_{i_k,j_l}$ of $\liea{gl}_N$ 
since every entry $(k,l)$ of such a rectangle 
corresponds uniquely to a matrix 
entry $E_{i_k,j_l}$. 
 
\begin{re} 
Our choice of the subset $S_1$ of the roots of $\liea{g}_1$ has 
to do with the weight structure of the representation 
of the Levi factor on the space $\liea{g}_1$. The subsets 
$S_1$ are related to our initial proofs of the theorems. 
Especially, they play a role in our first theoretical 
approach to the exceptional case. 
\end{re} 
 
\begin{re} 
Using Proposition~\ref{prop:mainBCD} one checks that the 
constructed nilpotent element is generic 
(that is the $X_R^j$ have generic rank). 
So the constructed matrix $X_R$ is a 
Richardson element for $\liea{p}$ in $\liea{g}_1$ 
\end{re} 
 
%
%
\subsection{Case $\lieg{A}_n$}\label{sse:A_n} 
%
%
Let $\liea{p}$ be a nice parabolic subalgebra 
of $\liea{sl}_{n+1}$. Let $(a_1,\dots, a_{r+1})$ 
be the sequence of the block lengths of the 
standard Levi factor of $\liea{p}$. 
By Theorem~\ref{th:A} this sequence is unimodal, so 
$a_1\le\dots\le a_s\ge\dots\ge a_{r+1}$. 
We will now give the subset $S_1$ of the sets of 
roots of $\liea{g}_1$ by explicitely describing 
the needed roots for each rectangle $R_{i,i+1}$. 
 
\begin{recipe}\label{recipe:A} 
(1a) If $i$ is odd we choose the entries $(1,1)$, 
$(2,2)$, $\dots, (a_i,a_i)$ of the rectangle 
$R_{i,i+1}$. 
 
(1b) If $i$ is even we choose entries 
$(a_{i+1},a_i)$, $(a_{i+1}-1,a_i-1)$, 
$\dots,(a_{i+1}-a_i+1,1)$ of $R_{i,i+1}$. 
 
(2) The subset $S_1$ of the roots of $\liea{g}_1$ 
is then the union of all roots corresponding 
to the chosen entries. 
 
(3) 
We define the nilpotent element $X_R$ to be the 
matrix $\sum_{\alpha\in S_1}E_{\alpha}$. 
\end{recipe} 
 
The result is a nilpotent matrix 
that has $r$ small boxes, one 
in each of the rectangles $R_{i,i+1}$. 
The small boxes consist of skew-diagonal matrices 
as shown in the following example.

\begin{ex}\label{ex:A_6} 
The parabolic subalgebra $(1,0,1,0,0,1)$ 
of $\lieg{A}_6$, i.e. 
$\liea{p}$ is given by the sequence 
$(1,2,3,1)$ of block length of the standard 
Levi factor. The normal form of the Richardson 
element is 
 
$X_R:=$ 
\begin{scriptsize} 
$\begin{pmatrix} 
* & 1 &  \\ 
 & * & * & & &1 \\ 
 & * & * &  & 1&\\ 
 & & & * & * & * &  \\ 
 & & & * & * & * &  \\ 
 & & & * & * & * & 1  \\ 
 & & & & & & * 
\end{pmatrix}$. 
\end{scriptsize} 
 
The roots that are chosen are 
$\alpha_1,\alpha_2+\alpha_3+\alpha_4,\alpha_6$, 
$\alpha_3+\alpha_4$. 
\end{ex} 
\begin{re} 
In Example 4.6 of~\cite{gr} R\"ohrle and Goodwin 
describe Richardson elements of $\liea{g}_1$ for parabolic 
subalgebras of $\liea{gl}_N$. They follow a construction 
given by~\cite{bhrr}. The method amounts to choose an 
identity matrix of size $\min\{a_i,a_{i+1}\}$ in each rectangle 
$R_{i,i+1}$, starting with entry $(1,a_i)$. The two methods are 
equivalent. Our approach is motivated by the construction for the 
other classical Lie algebras as we will see. 
\end{re} 
%
%
%
\subsection{Case $\lieg{C}_n$}\label{sse:Cn} 
%
 
Here, the pattern is slightly different. 
Either there is an even number of blocks in the 
standard Levi factor ((A)-case) or an odd number 
((B)-case). 
In the first case the recipe is essentially 
the same as for $\liea{sl}_{n+1}$. 
In the second case 
we modify Recipe~\ref{recipe:A}. 
We divide the chosen entries into two sets, 
choosing half of them starting with entry $(1,1)$ 
and the other half starting with entry 
$(a_{i+1},a_i)$. 
Let $b_i:=\lfloor\frac{a_i}{2}\rfloor$ 
and $B_i:=\lceil\frac{a_i}{2}\rceil$. 
 
\begin{recipe}\label{recipe:C} 
(1A) For $i\le r$ we do the same as in 
Recipe~\ref{recipe:A}. 
In particular, in rectangle $R_{r,r+1}$ we choose 
all entries on the skew-diagonal, i.e. 
$(1,1)$, $\dots,(a_r,a_r)$. 
 
(1B) 
If $i$ is odd, we choose the 
entries $(1,1)$, $\dots,(B_i,B_i)$ 
and the entries $(a_{i+1},a_i)$, 
$\dots,(a_{i+1}-b_i+1,a_i-b_i+1)$ 
of the rectangle $R_{i,i+1}$. 
If $i$ is even, we choose the 
entries $(1,1)$, $\dots,(b_i,b_i)$ 
and the entries $(a_{i+1},a_i)$, 
$\dots, (a_{i+1}-b_i,a_i-b_i)$. 
 
(2) 
The subset $S_1$ of the roots of $\liea{g}_1$ 
is the union of all roots that correspond 
to the chosen entries. 
 
(3) 
We define the nilpotent element $X_R$ to 
be the following matrix: 
If $\alpha_{ij}$ (with $i\le n$) is an element 
of $S_1$, then 
the matrix $X_R$ has entry $1$ at 
position $(i,j)$ and the correspond entry at its 
adjoint position: 
It has entry $-1$ at position 
$(2n-j+1,2n-i+1)$ if $j>n$ and entry $1$ at 
position $(2n-j+1,2n-i+1)$ if $j\le n$. 
\end{recipe} 
 
It is clear that this is a matrix of $\liea{sp}_{2n}$. 
 
\begin{ex} 
\begin{enumerate} 
\item 
The parabolic subalgebra $(0,0,1,0,0)$ 
of $\lieg{C}_5$. The standard Levi factor 
has block lengths $(3,4,3)$. There is an 
odd number of blocks. 
 
$X_R=$ 
\begin{scriptsize} 
$\begin{pmatrix} 
*&*&*&&&&1\\ 
*&*&*&&1&\\ 
*&*&*&1&&&\\ 
&&&*&*&*&*& &&\phantom{-}1\\ 
&&&*&*&*&*& & \\ 
&&&*&*&*&*& & -1\\ 
&&&*&*&*&*&-1 \\ 
&&&&&&&*&*&* \\ 
&&&&&&&*&*&* \\ 
&&&&&&&*&*&* 
\end{pmatrix}$ 
\end{scriptsize}. 
 
We have 
$\dim\ker X_R=4$, $\dim\ker X_R^2=8$. Thus the 
partition of $X_R$ is $(2,2,3,3)$, with one pair of 
odd entries. The dual partition is $(4,4,2)$. This gives 
$\dim\liea{g}^{X_R}=\frac{1}{2}(16+16+4)+1=19$ 
(cf. Theorem~\ref{thm:centralizer-dim-C} which 
is equal to $\dim\liea{m}$. 
 
\item 
The parabolic subalgebra $(0,1,0,0,0,1)$ 
of $\lieg{C}_6$. The standard Levi factor 
has block lengths $(2,4,4,2)$, i.e. an even number 
of blocks. 
 
$X_R=$ 
\begin{scriptsize} 
$\left( 
\begin{array}{cccccccccccc} 
*&*&&1\\ 
*&*&1\\ 
&&*&*&*&*& & & & 1\\ 
&&*&*&*&*& & &1\\ 
&&*&*&*&*& &1 \\ 
&&*&*&*&*&1 \\ 
&&&&&&*&*&*&*\\ 
&&&&&&*&*&*&*  \\ 
&&&&&&*&*&*&*& &-1 \\ 
&&&&&&*&*&*&*&-1 \\ 
&&&&&&&&&&*&* \\ 
&&&&&&&&&&*&* 
\end{array} 
\right)$ 
\end{scriptsize}. 
 
We compute the dimension of the centralizer 
of $X_R$ in $\liea{g}$. 
Note that $\dim\ker X_R=12-8=4$, $\dim\ker X_R^2=8$, 
$\dim\ker X_R^3=10$ and $\dim\ker X_R^4=12$. 
Thus the partition is $(2,2,4,4)$, which is 
equal to its dual partition. So 
$\dim\liea{g}^{X_R}=\frac{1}{2}(16+16+4+4)=20$. 
The Levi factor also has dimension $20$. 
\end{enumerate} 
\end{ex}

%
\subsection{Cases $\lieg{B}_n$, $\lieg{D}_n$} 
%
%
The choice in the case of the orthogonal Lie 
algebras is similar to the symplectic 
case. There are some differences: 
 
\begin{re}\label{re:adaptionsBD} 
In the construction for the special orthogonal 
case, there are some adaptions 
of the recipe: 
If there is an even number of blocks, i.e. in the 
(A)-case, and $a_i>a_{i+1}$ for some $i$, we 
pick $a_{i+1}$ roots for $R_{i,i+1}$. 
In the (B)-case we choose twice 
$B_i:=\lceil\frac{a_i}{2}\rceil$ entries 
in each rectangle $R_{i,i+1}$ if the parities 
of $a_i$ and of $a_{i+1}$ are different. 
\end{re} 
 
Recall that $b_i:=\lfloor\frac{a_i}{2}\rfloor$ 
and $B_i:=\lceil\frac{a_i}{2}\rceil$. To deal 
with non-unimodal sequences 
we introduce $s_i:=\min\{a_i,a_{i+1}\}$. 
 
\begin{recipe}\label{recipe:BD} 
(1A) 
Let the number of blocks in the standard 
Levi factor be even. 
 
Let $i\le r+1$. 
For odd $i$ 
we choose entries 
$(1,1)$, $\dots,$ $(s_i,s_i)$ in 
rectangle $R_{i,i+1}$. 
For even $i$ 
we choose entries 
$(a_{i+1},a_i),$ $\dots,$ 
$(a_{i+1}-s_i+1,a_i-s_i+1)$ of the 
rectangle $R_{i,i+1}$. 
 
In rectangle $R_{r,r+1}$ we choose 
$b_r$ two-by-two matrices 
$\begin{scriptsize} 
\begin{pmatrix} 1 & 0 \\ 0&-1\end{pmatrix} 
\end{scriptsize}$ on the skew-diagonal, 
starting at $(1,1)$ 
if $r$ is odd or at $(a_r,a_r)$ if $r$ is even. 
Cf. Example~\ref{ex:D-even}. 
 
(1B) 
For each rectangle $R_{i,i+1}$ we pick the entries 
$(1,1),\dots,(B_i,B_i)$ and the entries 
$(a_{i+1},a_i)$, $\dots,$ 
$(a_{i+1}-b_i,a_i-b_i)$. 
As an illustration, see 
Example~\ref{ex:BD-unim} below and 
Example~\ref{ex:BD-nonu}. 
 
(2) 
The set $S_1$ is the union of all roots 
$a_{ij}$ that correspond to the chosen 
entries. 
 
(3) 
The nilpotent element $X_R$ is the following 
matrix: 
If $\alpha_{ij}\in S_1$ then $X_R$ has entry 
$1$ at position $(i,j)$ and entry $-1$ at 
the position adjoint to $(i,j)$: at position 
$(2n-j+2,2n-i+2)$ in the case of $\liea{so}_{2n+1}$ 
resp. at $(2n-j+1,2n-i+1)$ in the case of 
$\liea{so}_{2n}$. 
\end{recipe} 
 
Example~\ref{ex:BD-nonu} treats 
non-unimodal cases. 
 
\begin{ex}\label{ex:BD-unim} 
We present the constructed nilpotent element for 
 parabolic subalgra of $\lieg{B}_4$ 
where the sequence of block lengths is unimodal 
and consists of an odd number of blocks. 
The parabolic subalgebra 
of $\lieg{B}_4$ is $(0,1,0,0)$, with block 
lengths $(2,5,2)$. 
 
$X_R=$ 
\begin{scriptsize} 
$\left( 
\begin{array}{ccccccccccc} 
*&*&*&&&&&1\\ 
*&*&*&&1&&1\\ 
*&*&*&1&&&&\\ 
&&&*&*&*&*&*& & &-1 \\ 
&&&*&*&*&*&*& &-1\\ 
&&&*&*&*&*&*& &\\ 
&&&*&*&*&*&*& &-1 \\ 
&&&*&*&*&*&*& -1 \\ 
&&&&&&&&*&*&* \\ 
&&&&&&&&*&*&* \\ 
&&&&&&&&*&*&* 
\end{array}\right)$ 
\end{scriptsize} 
 
Here, $\dim\ker X_R=5$, $\dim\ker X_R^2=8$ and 
$\dim\ker X_R^3=11$, giving the partition $(3,3,3,1,1)$ 
with five odd parts. 
The dual partition is $(5,3,3)$, and 
so 
$\dim\liea{g}^{X_R}=\frac{1}{2}(25+9+9)-\frac{5}{2}$ 
$=19$ which is the dimension of the Levi factor. 
\end{ex} 

\begin{ex}\label{ex:D-even} 
This matrix is the constructed nilpotent 
element for the parabolic subalgebra 
of $\lieg{D}_5$ 
that has $\alpha_2$ and $\alpha_5$ in 
the nilpotent radical. The 
sequence of the block lengths is 
$(2,3,3,2)$ with an even number of blocks. 
 
$X_R=$ 
\begin{scriptsize} 
$\left( 
\begin{array}{cccccccccc} 
*&*&&1&\\ 
*&*&1&&\\ 
&&*&*&* && 1&\phantom{-}0 \\ 
&&*&*&* & & 0&-1\\ 
&&*&*&* & & & \\ 
&&&&&*&*&\phantom{-}*&\\ 
&&&&&*&*&\phantom{-}* & &-1  \\ 
&&&&&*&*&\phantom{-}*&-1 \\ 
&&&&&&&&*&* \\ 
&&&&&&&&*&* 
\end{array} 
\right)$ 
\end{scriptsize} 
 
We have $\dim\ker X_R=4$, $\dim\ker X_R^2=6$ 
and $\dim\ker X_R^3=8$, $\dim\ker X_R^4=10$. 
The Jordan form is 
given by the partition $(1,1,4,4)$ with dual 
partition $(4,2,2,2)$. So the centralizer of 
$X_R$ has dimension $13$ which is 
equal to $\dim\liea{m}$. 
\end{ex} 
\begin{ex}\label{ex:BD-nonu} 
The parabolic subalgebra $(0,0,0,1,0)$ 
of $\lieg{B}_5$. The sequence of block 
length in the standard Levi factor is 
$(5,3,5)$. 
 
$X_R=$ 
\begin{scriptsize} 
$\left(\begin{array}{ccccccccccc} 
*&*&*&*& &&1\\ 
*&*&*&*& &1\\ 
*&*&*&*& &1\\ 
*&*&*&*&1&&&&\\ 
&&&&*&*&*& & & & -1\\ 
&&&&*&*&*& &-1&-1 \\ 
&&&&*&*&*&-1& & \\ 
&&&&&&&*&*&*&*\\ 
&&&&&&&*&*&*&* \\ 
&&&&&&&*&*&*&* \\ 
&&&&&&&*&*&*&* 
\end{array}\right)$ 
\end{scriptsize} 
 
We leave the calculation of the dimension of 
the Levi factor to the reader. 
The second example is the parabolic subalgebra 
of $\lieg{D}_6$ where $\alpha_2$, $\alpha_3$ 
and $\alpha_6$ 
are the simple roots in the Levi factor. 
The sequence of the block lengths is 
$(1,3,2,2,3,1)$, with an even number of blocks. 
 
$X_R=$ 
\begin{scriptsize} 
$\left(\begin{array}{cccccccccccc} 
*&1 \\ 
&*&*&*& &1 \\ 
&*&*&*&1 \\ 
&*&*&*& \\ 
& & & &*&*&1& \\ 
& & & &*&*& &-1 \\ 
& & & & & &*&*& & &-1\\ 
& & & & & &*&*& &-1 \\ 
& & & & & & & &*&*&*& \\ 
& & & & & & & &*&*&*& \\ 
& & & & & & & &*&*&*&-1\\ 
& & & & & & & & & & &* 
\end{array}\right)$ 
\end{scriptsize} 
\end{ex} 
 
\begin{re} 
In the (A)-cases, type $\lieg{B,D}$ Recipe~\ref{recipe:BD} 
chooses more roots than we do in 
the other case. 
We know that there are example where this is 
necessary, e.g. 
the parabolic subalgebra of $\lieg{D}_{11}$ 
with block lengths $(1,3,5,4,5,3,1)$. 
Here, one needs the ``extra'' 
roots chosen in $R_{1,2}$ and in $R_{2,3}$. 
If we choose only the commuting roots according to 
the recipe for type $\lieg{C}$ (Recipe~\ref{recipe:C}) 
we obtain 
a nilpotent $X_R'$ with $\dim\liea{g}^{X_R'}>41$ 
while $\dim\liea{m}=41$. 
On the other hand, if there are only odd sized 
blocks, as in Example~\ref{ex:BD-unim}, there is 
a way to choose one entry less per rectangle 
$R_{i,i+1}$. In Example~\ref{ex:BD-unim} 
one can alternatively choose $(2,1)$, $(3,2)$ 
and $(4,3)$ in $R_{1,2}$. 
Also, if $r=2$, we can use part (1B) of Recipe~\ref{recipe:C} instead of 
(1B) in Recipe~\ref{recipe:BD}. 
These phenomena are discussed in detail in~\cite{ba} where 
the first named author describes 
the properties of the roots involved in the construction 
of Richardson elements. 
\end{re}

\end{document}